\documentclass[10pt,a4paper,twoside,notitlepage]{article}

\usepackage{amsfonts}
\usepackage{amsmath}
\usepackage{amssymb}

\usepackage{graphicx}
\usepackage{psfrag}
\usepackage{color}

\newtheorem{thm}{Theorem} 
\newtheorem{lem}{Lemma}[section]
\newtheorem{cor}{Corollary}[section]
\newtheorem{prop}{Proposition}[section]

\newtheorem{rem}{Remark}[section]

\newcommand{\R}{\mathbb{R}}
\newcommand{\N}{\mathbb{N}}

\newcommand{\ve}{\varepsilon}
\newcommand{\n}{\noindent}
\newcommand{\vp}{\varphi}

\newcommand{\Div}{\mathrm{div}}

\newcommand{\D}{\mathrm{D}}

\newcommand{\con}{\mathbf{C}}
\newcommand{\Lsp}{\mathbf{L}}

\newcommand{\meas}{\mathrm{meas}}

\title{Well-posedness of 2-$D$ and 3-$D$ swimming models in incompressible fluids governed by Navier--Stokes equations\footnote{This work  was supported in part  by the INdAM national group GNAMPA and the INdAM-CNRS GDRE CONEDP, and by  the Grant 317297 from Simons Foundation. }} 

\date{}

\begin{document}

\maketitle

\begin{center} 
Corresponding author:  Alexandre Khapalov\footnote{The work by this author was supported in part by NSF Grant DMS-1007981.} ,
Department of Mathematics, \\ 
Washington State University,  USA (email: khapala@math.wsu.edu)\\
Piermarco Cannarsa, Department of Mathematics, University of Rome ``Tor Vergata'',   Italy   \\
Fabio S. Priuli, Istituto per le Applicazioni del Calcolo ``M. Picone'' (C.N.R.), Rome, Italy \\
Giuseppe Floridia, Department of Mathematics, University of Rome ``Tor Vergata'',   Italy 
\end{center}

\begin{abstract}
We introduce and investigate the wellposedness of two models describing the self-propelled  motion of a ``small bio-mimetic swimmer'' 
in the  2-$D$ and 3-$D$ incompressible fluids modeled  
by the  Navier-Stokes equations.  It is  assumed that the swimmer's body  consists of 
finitely many subsequently connected parts, identified with the fluid they occupy,  linked by
 the rotational  and elastic  forces. The swimmer employs the change of its  shape, inflicted by respective explicit internal forces, as the means for self-propulsion in a surrounding medium.  Similar models were  previously investigated  in  \cite{Kh11}-\cite{Kh5}  where the fluid was modeled by  the liner nonstationary Stokes equations.
Such models are of interest in biological and engineering applications dealing with the study and design of propulsion systems in fluids and air. 
\end{abstract}


\bigskip


\section{Introduction}\label{sec:intro}
The swimming phenomenon has been the subject of  interest for many researchers in various areas of natural sciences for a long time,  aimed primarily at understanding  biomechanics of swimming locomotion of biological organisms, see Gray
 \cite{Gra1}(1932), Gray and Hancock  \cite{Gra2} (1951), Taylor  \cite{Tay1} (1951),  \cite{Tay2} (1952),  Wu \cite{Wu} (1971), Lighthill  \cite{Lig} (1975),
and others. This research  resulted in the
derivation of a number of mathematical models for swimming motion  in 
the (whole)  ${\mathbb R^2}$- or ${\mathbb R^3}$-spaces with the  swimmer to be used as 
the reference frame, see, e.g.,
Childress  \cite{Chi}  (1981) and the references therein. In particular,  based on the
size of Reynolds  number, it was suggested (for the purpose of simplification)  to divide swimming models into three groups: microswimmers  (such as
flagella, spermatozoa,  etc.) with ``insignificant'' inertia;
``regular'' swimmers (fish,
dolphins, humans, etc.), whose motion takes
into account both  viscosity of  fluid and inertia; and Euler's swimmers,\index{Euler's swimmer} in which case
viscosity is to be `'neglected''.

\n
It also appears that the following two,  in fact, mutually excluding approaches were distinguished to model the swimming phenomenon (see, e.g., Childress, \cite{Chi}  (1981)). One,  which we can call the ``{\em shape-transformation approach}'',  exploits the idea that the swimmer's shape transformations during the actual swimming process  can be viewed as  a set-valued map in time (see  the seminal paper by  Shapere and Wilczeck  \cite{Sha} (1989)). The respective models describe the  swimmer's position in a fluid via the aforementioned  maps, see, e.g.,   \cite{Cur} (1981),   \cite{San2} (2008), \cite{Dal} (2011)  and the references therein.  Typically, such  models consider these maps as {\em a priori prescribed}, in which case the question  whether the respective  maps  are admissible, i.e., compatible with the principle of self-propulsion  of swimming locomotion or not, remains unanswered. In other words, one cannot guarantee that the model at hand describes the respective motion as a self-propulsive, i.e.,  swimming process. 
To ensure the positive answer to this question one needs to be able to answer the question whether the a priori prescribed  body changes of swimmer's shape can indeed be a result of actions of its internal forces under unknown in advance interaction with the surrounding medium.

\n
The other modeling approach (we will call it the {\em ``swimmer's internal forces approach'' or SIF approach}) assumes that the available internal swimmer's forces are explicitly described in the model equations and, therefore,  they determine the  resulting swimming motion. In particular, these forces will define the respective   swimmer's shape transformations in time as a result of an unknown-in-advance interaction of swimmer's body with the surrounding medium under the action of the aforementioned forces. For this approach,  we refer to   Peskin  \cite{Pes1}  (1975), Fauci and Peskin  \cite{Fau1} (1988), Fauci  \cite{Fau2}  (1993), Peskin and McQueen  \cite{Pes2}  (1994), Tytell, Fauci et al \cite{Fau3}  (2010), Khapalov \cite{KhBook}  and the references therein.

\n
The original idea of Peskin's approach is to view a  ``narrow'' swimmer as  an immaterial  ``immersed boundary''.  Within this approach the swimming motion is defined at each moment of time by the explicit swimmer's internal forces.
Following the ideas of Peskin's approach, Khapalov introduced the {\em immersed body SIF  modeling approach}   in which the bodies of ``small'' flexible swimmers are assumed to be   identified with the fluid  within their shapes, see   \cite{Kh11}-\cite{Kh5} (2005-2014). Indeed, in the framework of Peskin's method the swimmer  is modeled as  an immaterial curve, identified with the fluid,  further discretized for computational purposes on some grid as a collection of finitely many ``cells'', which in turn can be viewed as an immerse body, see Figures 1 and 2.
The idea here is to try, making use of mathematical simplifications of such approach,  to   focus  on the issue of macro dynamics of a swimmer. The simplifications (they seem to us to be legitimate within the framework of our interest) include the  reduction of  the number of model equations  and avoiding  the analysis of micro level interaction between a ``solid'' swimmer's body surface and fluid. It should be noted along these lines that in typical  swimming models dealing with ``solid'' swimmers, the latter are modeled  as ``traveling holes'' in system's space domain, that is, the aforementioned ``micro level'' surface interaction is  not in the picture as well.  

\n
In the above-cited works  by Khapalov \cite{Kh11}-\cite{Kh5},  the immerse body SIF approach was applied to the nonstationary Stokes equations in 2-$D$ and 3-$D$ dimensions with the goal to investigate the well-posedness of  respective models and their controllability properties.
In this paper our goal is to extend these results with respect to well-posedness  to  the case of Navier-Stokes equations in both the 2-$D$ and 3-$D$ incompressible fluids. To our knowledge, there were no previous publications investigating this issue within the SIF  approach.

\begin{figure}[h]\label{fig:swim2d}
\centering
\includegraphics[scale=0.6]{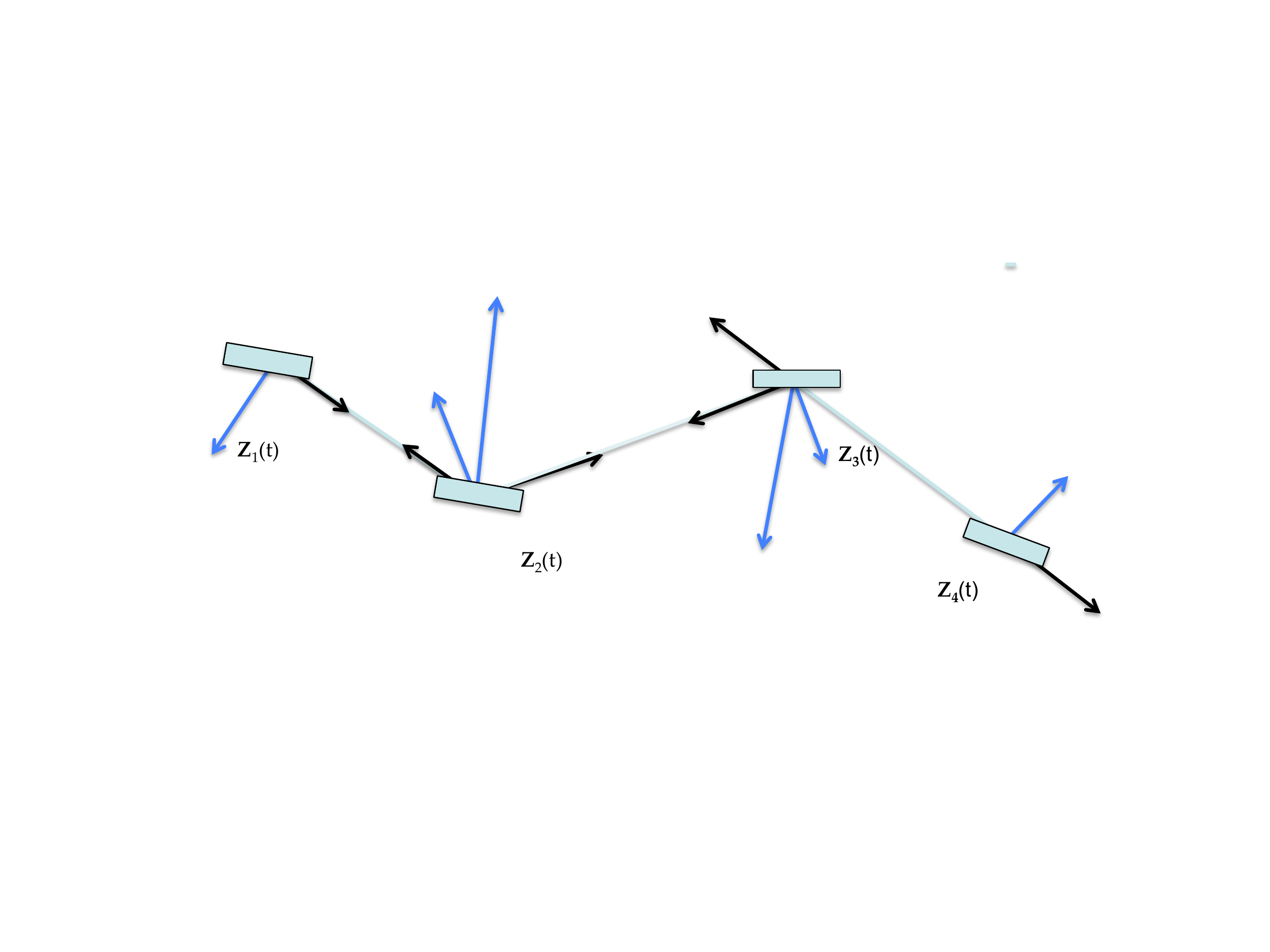} 
\caption{2-D swimmer consisting of 4 rectangles.}
\end{figure}

\begin{figure}[h]\label{fig:swim3d}
\centering
\includegraphics[scale=0.6]{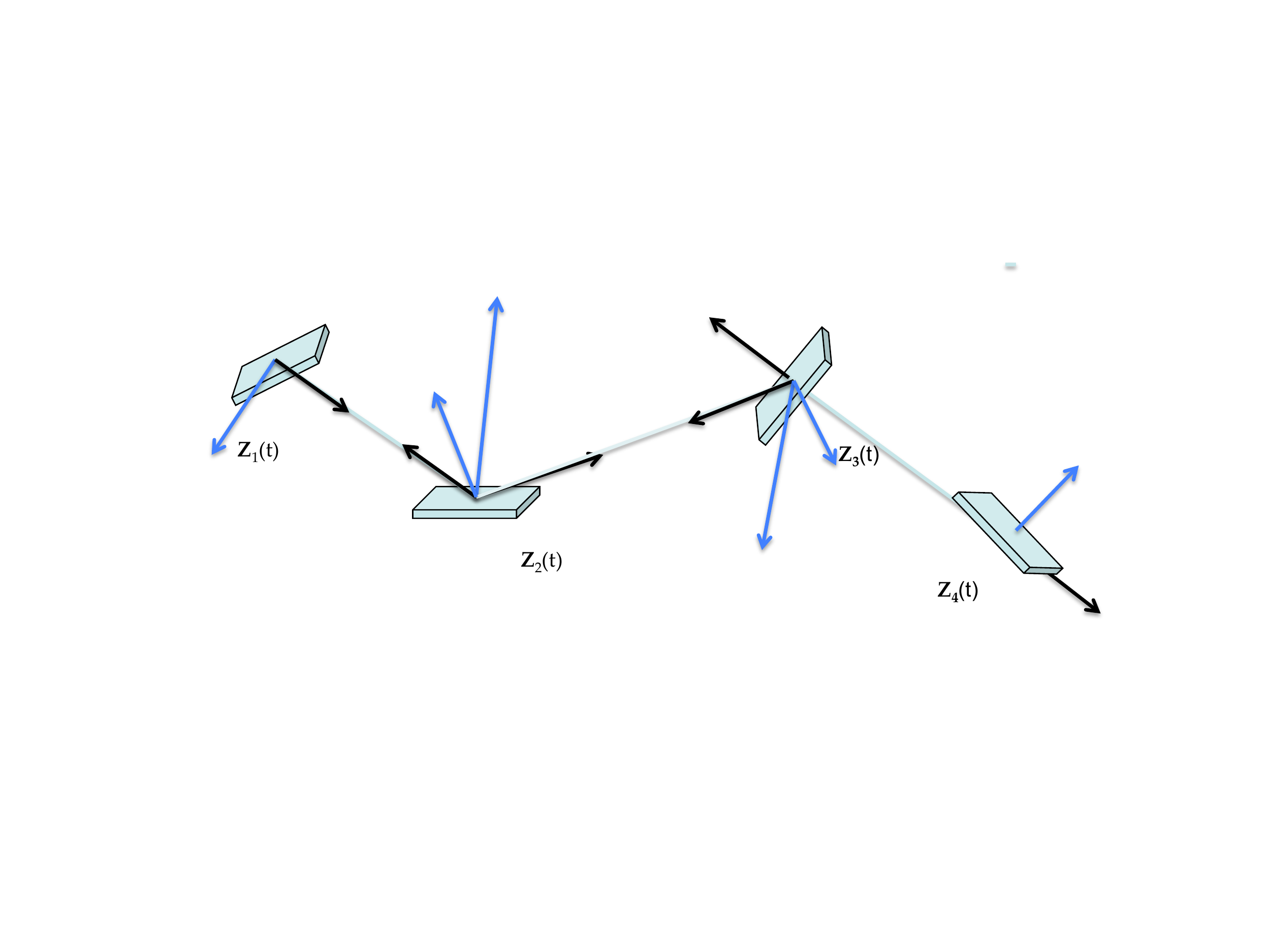} 
\caption{3-D swimmer consisting of 4 parallelepipeds.}
\end{figure}

\n
{\bf Related references on well-posedness of swimming models.}  To our knowledge,  in the context of PDE approach to swimming modeling,  the classical mathematical issues of  well-posedness   were addressed   for the first time  by
Galdi  \cite{Gal}  (1999), \cite{Gal2}  (2002) for a model of swimming micromotions  in $ {\mathbb R^3}$
(with the swimmer as the reference frame).
In \cite{San2} (2008)  San Martin et al discussed the well-posedness of  a 2-$D$ swimming model  within the framework of the shape transformation approach for the fluid governed by the 2-$D$  Navier-Stokes equations.

\n
{\bf Swimming models in the framework of ODE's.} A number of  attempts were made to introduce various reduction techniques  to convert swimming model equations into systems of ODE's (e.g.,  by making use of  empiric  observations and experimental data, etc.), see, e.g., 
Becker et al  \cite{Bec}  (2003); Kanso et al. \cite{Kan} (2005);  Alouges et al.  \cite{Alo} (2008), Dal Maso et al. \cite{Dal} (2011) and the references therein.

\bigskip

\n
The paper is organized as follows. In Section 2 we state our main results. In Section 3 we discuss in detail the modeling approach of this paper to swimming locomotion. In Section 4 we prove our main results after stating  several auxiliary lemmas, proven further in Section 5.
In Sections 6 and 7 (Appendices A and B) we remind the reader some classical results used in our proofs.


\section{Main results}\label{sec:main}

To formulate these results, we will need the following function spaces.

\bigskip
\n
Let $\Omega\subseteq\R^d$ be a bounded domain with locally Lipschitz boundary $\partial \Omega$. Below,  we use the following classical notations:
\begin{itemize}
\item $d$ denotes the dimension of the space domain,  equal either to 2 or to 3;
\item $\con_c^\infty(\Omega)$ denotes  the space of infinitely many times differentiable functions with compact support  in $\Omega$;
\item ${\cal D}'(\Omega)$ denotes  space of distributions in $\Omega$, i.e., the dual space of $\con_c^\infty(\Omega)$;
\item $W^{1,p}(\Omega)$, $1\leq p\leq \infty$  denote the Sobolev spaces over $\Omega$, i.e., the Banach spaces of functions in $\Lsp^p(\Omega)$ whose first (generalized) derivatives belong to $\Lsp^p(\Omega)$;
\item $H^1(\Omega)=W^{1,2}(\Omega)$, and  $H^2(\Omega)=W^{2,2}(\Omega)$;
\item $H^1_0(\Omega)$ denotes the subspace of   $H^1(\Omega)$ consisting of functions vanishing on $ \partial \Omega$.
 $H^{-1}(\Omega)$ denotes the dual space of $H^1_0(\Omega)$.
\end{itemize} 
Following~\cite{Temam}, page 5, we also introduce the following  $d$-dimensional vector function spaces:
$$
{\cal V}:=\{\vp\in[\con_c^\infty(\Omega]^d~;~\Div\,\vp=0\}\,,
$$
$$
H:=\mathrm{cl}_{\Lsp^2}({\cal V})\,,
\qquad\qquad
V:=\mathrm{cl}_{H^1_0}({\cal V})\, =\{\vp\in [H^1_0(\Omega)]^d~;~\Div\,\vp = 0\}\,,
$$
where the symbol $\mathrm{cl}_{\Lsp^2}$ stands for  the closure with respect to the $[\Lsp^2(\Omega)]^d$-norm, and $\mathrm{cl}_{H^1_0}$ -- with respect to the $[H^1_0(\Omega)]^d$-norm. The latter is  induced  by the scalar product
$$
\ll\!\vp,\psi\!\gg~:=~\sum_{j=1}^d \langle\D_j\vp, D_j\psi\rangle_{\Lsp^2}\ = \; \sum_{i,j=1}^d \mathop{\int}_\Omega \vp_{i x_j} \psi_{i x_j} dx,
\; \; \vp = (\vp_1,  \ldots, \vp_d), \; \psi = (\vp_1,\ldots,  \vp_d),
$$
where $ D_j$ is the differentiation operator with respect to $ x_j$.
To simplify notations, below we will use the notation $\|\vp\|_{\Lsp^2}$ (resp. $\|\vp\|_{H^1_0}$) both for functions $\vp\in \Lsp^2(\Omega)$ (resp.  $\vp\in H^1_0(\Omega)$) and for functions $\vp\in [\Lsp^2(\Omega)]^d$ (resp.  $\vp\in [H^1_0(\Omega)]^d$).

\bigskip

\n
Let $V'$ and $H'$ stand for  the dual spaces  respectively of $V$ and $H$. Then,  identifying $ H$ with $ H^\prime$, we have
$$
V\subset H\equiv H'\subset V'\,.
$$

\medskip
\n
Our main results deal with the well-posedness of 2-$D$ and 3-$D$ swimming models (\ref{eq:nse}), (\ref{eq:ode}), (\ref{eq:forces})-(\ref{eq:forces_rot3}),  described in detail in the next section and visualized by  Figures 1 and 2.

\begin{thm}[Well-posedness of the 2-$D$ swimming model]\label{thm:wellposed}
Let $d=2$ and for some $ T~>~0$  let $ \; u_0 \in H$,  $z_{1,0},\ldots,z_{N,0} \in \Omega$, $\kappa_1,\ldots,\kappa_{N-1},v_1,\ldots,v_{N-2} \in \Lsp^2(0, T)$. Assume that the assumptions {\bf (H1)}-{\bf (H2)} (given in Section 3) hold, and that
 \begin{equation}\label{body1}
\overline{S}(z_{i,0}) \subset \Omega, \qquad\qquad
|z_{i,0}-z_{j,0}| > 2r\,, \qquad\qquad i, j = 1, \ldots, N,~~ i \neq j\,,
\end{equation} 
where $r>0$ is the constant in {\bf (H1)}. (Assumption (1) ensures that no parts of swimmer's body overlap with each other and all lie within $\Omega$.)
Then, there exists $T^*\in(0,T]$ such that system~(\ref{eq:nse}), (\ref{eq:ode}), (\ref{eq:forces}),(\ref{eq:forces_el})/(\ref{eq:forces_el2}),(\ref{eq:forces_rot2})  admits a unique solution $(u, z)$ in $C(0,T^*\,;H)\cap\Lsp^2(0,T^*\,;V)\times [\con([0,T^*]\,;\R^d)]^N$, and 
 \begin{equation}\label{body2}
\overline{S}(z_i(t)) \subset \Omega, \qquad
|z_i(t)-z_j(t)| > 2r\,, \qquad i, j = 1, \ldots, N,~~ i \neq j \, \qquad \forall  t \in [0, T^*].
\end{equation}
The formula for  $\nabla p$, complementing the given pair $(u,z)$ in models (\ref{eq:nse}), (\ref{eq:ode}), (\ref{eq:forces})-(\ref{eq:forces_rot3}), is given   in Proposition \ref{prop:pressure} below.
\end{thm}

\begin{thm}[Well-posedness of the 3-$D$ swimming model]\label{thm:wellposed3D1} Let  $d=3$, $\partial \Omega$ be of class $ C^2$ and $u_0 \in  V$. 
Then the  result stated in Theorem 1 holds for model (\ref{eq:nse}), (\ref{eq:ode}), (\ref{eq:forces}),(\ref{eq:forces_el})/(\ref{eq:forces_el2}), (\ref{eq:forces_rot3})
for a unique triplet   $ (u, p, z) $ such that $ u_t, \Delta u, \nabla p \in [\Lsp^2 (Q_{T^*})]^3$ and $ u \in C([0,T^*]; V)$, where $Q_{T^*} = (0,T^*) \times \Omega$.  
\end{thm}

\begin{thm}[Additonal regularity]\label{thm:maddreg}
Let $\partial \Omega$ be of class $ C^2$. If  $u_0 \in [H^{2} (\Omega)]^3\bigcap V$, then  in Theorem \ref{thm:wellposed3D1} solution   $ u $ lies in $  [H^{2,1} (Q_{T^*})]^3 \bigcap C([0,T^*]; V)$.   In turn, if  $u_0 \in [H^{2} (\Omega)]^2\bigcap V$, then   $ u \in [H^{2,1} (Q_{T^*})]^2 \bigcap C([0,T^*]; V)$ in Theorem \ref{thm:wellposed}.
\end{thm}
Here, $ H^{2,1} (Q_{T^*}) = \{ \vp   \; | \; \vp  \in L^2 (0,T^*; H^2 (\Omega)), \; \vp_t \in L^2 (Q_{T^*})\}$.
This result is an immediate consequence of Theorem \ref{thm:addreg} in Section 6 (Appendix A).

\begin{rem}
The duration of time $T^*$ (i.e., of existence of solutions)  in Theorems \ref{thm:wellposed}-\ref{thm:maddreg} depends on the parameters of respective model $ u_0, z_{i,0}$'s, $k_i$'s, $v_i$'s and the initial shape and position of the swimmer in $ \Omega$.  One can view $ T^*$ as a new initial instant of time to apply these theorems again to further extend the interval of existence from $ T^*$ to some $T^{**} > T^*$, provided that the assumptions of the corresponding theorem hold at $ t =T^*$, and so on. Our proofs below indicate that in the 2-$D$ case one can extend this time up to the collision of swimmer with the boundary of $ \Omega$ or up to the moment when swimmer's body parts will collide with each other, see conditions (\ref{body1})-(\ref{body2}). In turn, these circumstances depend on or  can be regulated by a suitable choice of swimmer's internal forces, i.e.,  functions $v_i$'s (and $k_i$'s when they can vary, see the next section). To the contrary, in the 3-$D$ case, such interval of existence of solution will also depend on $ u_0, z_{i,0}$'s, $k_i$'s, and $v_i$'s via an additional condition (\ref{3dT1}) in  Theorem \ref{thm:existence3d} (compare it to  Theorem \ref{thm:existence2d}).
\end{rem}


\section{Swimming model}\label{sec:model}
Following  \cite{KhBook}, \cite{Kh4}-\cite{Kh5}, we describe the locomotion of a swimmer in a fluid by a hybrid  nonlinear  system of two sets of pde/ode equations:
 \begin{equation}\label{eq:nse}
\left\{\begin{array}{ll}
u_t-\nu\Delta u+(u\cdot\nabla)u+\nabla p  =f & \mbox{ in } (0,T)\times\Omega,\\
\Div \,u = 0& \mbox{ in } (0,T)\times\Omega,\\
u=0 & \mbox{ in } (0,T)\times\partial\Omega,\\
u(0,\cdot)=u_0 & \mbox{ in } \Omega,
\end{array}\right.
\end{equation}

\begin{equation}\label{eq:ode}
{dz_i\over dt}\,=\,{1\over \meas(S (0))}\,\int_{S(z_i(t))} u(t,x)\,dx,
\qquad 
z_i(0) = z_{i,0},
\qquad
i=1,\ldots,N, \qquad t \in (0, T).
\end{equation}
System (\ref{eq:nse})  describes the evolution of incompressible fluid due to the Navier-Stokes equations under the influence of the forcing term $ f (t,x)$ representing the actions of swimmer. Here,  $(u (t,x) ,p (t,x) )$ are respectively  the velocity of the fluid and its pressure at point $ x$ at time $t$, and $\nu$ is the kinematic viscosity constant.  In turn, system (\ref{eq:ode}) describes the motion of the swimmer  in $\Omega$, whose flexible body consists of $ N$   sets $ S (z_i (t))$ within $\Omega$. These sets  are identified with the fluid within the space they occupy at time $t$ and are linked between themselves by the  rotational and  elastic forces as illustrated on Figures 1 and 2. The points $ z_i (t)$'s represent the centers of mass of the respective parts of swimmer's body. The instantaneous velocity of each part  is calculated as  the  average fluid velocity within it at time $t$. 

\bigskip
\n
We will now describe the assumptions on the parameters of model (\ref{eq:nse})-(\ref{eq:ode}) in  detail.

\subsection{Swimmer's body}\label{sec:body_model}
Below, for simplicity of notations, we will denote the sets $S(z_i(t)), i=1,\ldots,N$ also  as $S(z_i)$ or $S_i(z_i)$. Throughout the paper we assume the following two main assumptions: 

\begin{description}
\item{\bf (H1)} {\em All sets $S(z_i ), i=1,\ldots,N$ are obtained by  shifting 
the same  set $S(0)\subset\Omega$, i.e.,
\begin{equation}\label{S0}
S(z_i)=z_i+S(0),\qquad\qquad i=1,\ldots,N\,, 
\end{equation}
where $S(0)$ is open and lies in a ball $B_r(0)$ of radius  $r>0$, and its center of mass is  the origin. }
\end{description}

\begin{rem}\label{orient}
The  results of this paper will hold at no extra cost if we assume that  the swimmer at hand consists of different body parts $ S_i(0) \subset B_r(0)$, in which case one will need to    replace (\ref{S0}) with $  S(z_i)= z_i+S_i(0), \;  i=1,\ldots,N$ and  add respective normalizing coefficients in the expressions for swimmer's internal forces 
(\ref{eq:forces})-(\ref{eq:forces_rot3}) to ensure that they satisfy the 3rd Newton's Law. In particular, the swimmers on Figures 1 and 2 consist of identical sets each of which has its own orientation   in space.
\end{rem}

\bigskip
\n
We will need the following concepts from ~\cite[Section 3.11]{AFP} to formulate our second  assumption on the geometry of  swimmer  in this paper,  which will also be used in the respective proofs below.

\n
Given a bounded set $\Omega\in\R^n$ and a unit vector  $\nu\in\R^n$, denote by
$$
\pi_\nu:=\{\xi\in\R^n~;~\xi\cdot\nu=0\}\,,
\qquad\qquad
\Omega_\nu:= \left\{\xi\in\pi_\nu~;~\exists\, t\in\R\mbox{ s.t. } \xi + t\nu\in\Omega \right\}\,
$$
respectively the orthogonal hyperplane to $\nu$ and the orthogonal projection of $\Omega$ on $\pi_\nu$. Then, for every $y\in\Omega_\nu$ we will call the set
\begin{equation}\label{eq:section_set}
\Omega_\nu^y:=\left\{t\in\R~;~y+t\nu\,\in\Omega, \; y\in\Omega_\nu\right\}
\end{equation}
 the  \emph{section} of $\Omega$ corresponding to $y$.  Accordingly, given a function $\phi \colon \Omega\to\R$ and any $y\in\Omega_\nu$, we define the function $\phi_\nu^y (t)$, called  \underline{\emph{section} {\em of $\phi$ corresponding to $y$}}, as
\begin{equation}\label{eq:section_map}
\phi_\nu^y\colon \R \supset  \Omega_\nu^y \ni t \to \phi_\nu^y(t) \in \R, \qquad \phi_\nu^y(t):= \phi\left(y+t\nu\right)\,.
\end{equation}
\begin{description}\item{\bf (H2)} {\em There exist positive constants $h_0$ and ${\cal K}_S$ such that for any vector $h\in B_{h_0}(0)\setminus\{0\}$ we can find a vector $\eta=\eta\big({h}\big)$, $\;  \mid \eta \mid   = 1$ which satisfies 
\begin{equation}\label{eq:H2}
\meas(S_\Delta)_\eta^y=\int_{(S_\Delta)_\eta^y}dt\leq {\cal K}_S\,|h|
\qquad\qquad\forall~y\in\Omega_\eta\,,
\end{equation}
where $S_\Delta:=(h+S(0)) \,\Delta\,  S(0)$ is the symmetric difference between $S(0)$ and $h+S(0)$, i.e. $S_\Delta=\big((h+S(0))\cup S(0)\big)\setminus \big((h+S(0))\cap S(0)\big)$.}
\end{description}

\smallskip
\n
Assumption {\bf (H2)} means that size of the projection of  $S_\Delta$ on the hyperplane in $ R^d$ perpendicular to vector  $\eta$ changes uniformly Lipschitz  continuously relative to the magnitude of the shift $ h$ of the set $S(0)$  in the direction of $h$. This assumption is principally weaker than  the respective assumption  on the regularity of the shifts of $ S(0)$ in \cite{KhBook}, \cite{Kh4}-\cite{Kh5}, where $ \eta$ was always selected  to be $h$ (we will illustrate it in  Example 3.1 below). In the case when $ \eta = h$, it is easy to verify that {\bf (H2)} is satisfied for  discs and  rectangles in 2-$D$ and for balls and parallelepipeds in 3-$D$.

\bigskip
\begin{rem}
In this paper we assume that all swimmer's body parts  are  identical sets. One can choose these sets to be of  distinct  shapes and sizes, in which case, however,  the respective normalizing coefficients should be added to the forcing terms to ensure that all swimmer's forces  are to be its internal forces.
\end{rem}

\n We conclude this subsection with an example showing that there exist particular shapes of the set $S(0)$ which require the presence of $\eta=\eta(h)$ in {\bf (H2)} instead of the straightforward choice $\eta={h\over|h|}$.

\n
{\bf Example 3.1. }
Fix a constant $\kappa>0$. We claim that there exist sets $S(0)$  which satisfies  {\bf (H2)} for some $ \eta$, but   $\eta$ cannot be selected to be co-linear with $ h$.
\begin{figure}[h]\label{fig:triang_sasha}
\centering
\includegraphics[width=.6\textwidth]{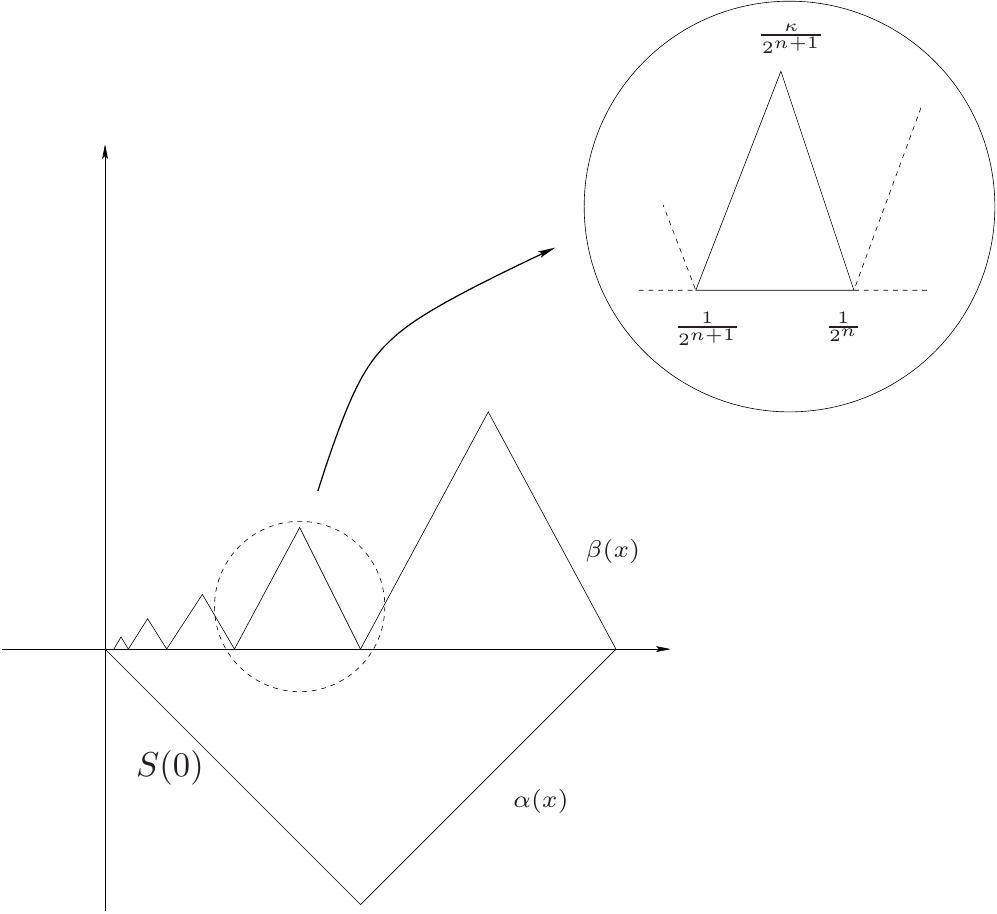} 
\caption{Set $S$ in  Example 3.1.}
\end{figure}

\n
Let us consider the  set $ S$ shown on Figure 3,
$$
S:=\left\{
(x,y)\in\R^2~|~ x\in (0,1)\,,\,y\in (\alpha(x),\beta(x)),
\right\}
$$
where the functions $\alpha,\beta\colon[0,1]\to\R$ are defined as follows:
$$
\alpha(x):=
\left\{
\begin{array}{ll}
-x, & \mbox{ if } \displaystyle 0 \leq x\leq\,{1\over 2},\\&\\
x-1, & \mbox{ if } \displaystyle 1 \geq x\geq\,{1\over 2},
\end{array}
\right.
\quad
\beta(x):=
\left\{
\begin{array}{ll}
\displaystyle 2\kappa \left(x-\,{1\over 2^{n+1}}\right), & \mbox{ if } \displaystyle  x\in\left[{1\over 2^{n+1}}\,,\,{3\over 2^{n+2}}\right],\\&\\
\displaystyle 2\kappa \left({1\over 2^n}\,-x\right), & \mbox{ if } \displaystyle  x\in\left[{3\over 2^{n+2}}\,,\,{1\over 2^n}\right],
\end{array}
\right.
$$
$n = 0, 1, \ldots$.
Observe that the part of the boundary $\partial S$ corresponding to $y=\beta(x)$ is given by lines with slopes either  $2\kappa$ or $-2\kappa$. If we introduce the notation ${\cal T}_m$ ($m\geq 1$) for the isosceles triangle of base $2^{-m}$ in $[2^{-m}, 2^{-m+1}]$ and height $\kappa/2^m$, and ${\cal T}_0$ for the isosceles triangle $\{(x,y)\in\R^2~;~ x\in [0,1]\,,\,y\in [\alpha(x),0]\}$, then we have
$$
\bar{S}=  \bigcup_{m\geq 0} {\cal T}_m\,.
$$
Denote by ${\frak b}(S)$ the center of mass of $S$ and set $S(0) = S-{\frak b}(S)$. 
We claim that, no matter what $h_0$ and ${\cal K}_S$ we choose, if we shift $S(0)$ by $h=\ve (-1,0)$, for a suitable $\ve\in(0,h_0)$, and if we use $\eta\big({h}\big)=(-1,0)$, then we can always find $\bar y\in \Omega_\eta$ such that
$$
\meas\, [S(0)~\Delta~(h+S(0))]_\eta^{\bar y} > {\cal K}_S\, \ve\,,
$$
and, thus, (\ref{eq:H2}) does not hold. On the other hand, by setting $\eta\big({h}\big)\equiv (0,1)$ for all $h\in B(0,1)$, we can prove that
\begin{equation}\label{eq:ex_H2_part1}
\meas\, [S(0)~\Delta~(h+S(0))]_\eta^y\leq 4 \kappa \ve \qquad\qquad \forall~y\in \Omega_\eta\,,
\end{equation}
so that ~\eqref{eq:H2} is satisfied with $h_0=1$ and ${\cal K}_S=4 \kappa$. 

\n
Since the Lebesgue measure is invariant with respect  to translations, in the computations below we will always use $S$ and  $S_\Delta=S~\Delta~(h+S)$ in place of $S(0)$ and $S(0)~\Delta~(h+S(0))$.

\smallskip
\n
We start with the negative result. For every fixed $h_0>0$ and ${\cal K}_S>0$, let us consider $m$ large enough to have $2^{-m}< h_0$ and $2(m-1)> {\cal K}_S $. By choosing $h= 2^{-m} (-1,0)\in B_{h_0}(0)$, we claim we can find $\bar y\in (S_\Delta)_h\subset\Omega_h$ such that $\meas(S_\Delta)_h^{\bar y}>{\cal K}_S\, |h|= {\cal K}_S/2^m$, so that the choice $\eta\big(h\big)={h\over|h|}$ is not suited for this set $S$.

\n 
The projection $(S_\Delta)_h$ on the $y$--axis (which is the orthogonal line to $h$) is the segment  $(-1/2, \kappa/2)$. Thus, the value $\bar y=2^{-m}\kappa$ belongs to $(S_\Delta)_h$ and (see Figures 3- 4)
$$
(S_\Delta)_h^{\bar y} = \bigcup_{1\leq k\leq m} \big({\cal T}_m~\Delta~(h+{\cal T}_m)\big)_h^{\bar y}\,.
$$
i.e.,  $ (S_\Delta)_h^{\bar y} $ consists of  the union of the sections of the symmetric differences of the triangles in $S$ whose height exceeds $2^{-m}\kappa$. Hence,  $\meas \, (S_\Delta)_h^{\bar y}\} >= 2(m-1) |h|>{\cal K}_S\,|h|$ from which our claim follows.

\begin{figure}[h]\label{fig:triang_sasha}
\centering
\includegraphics[scale=0.6]{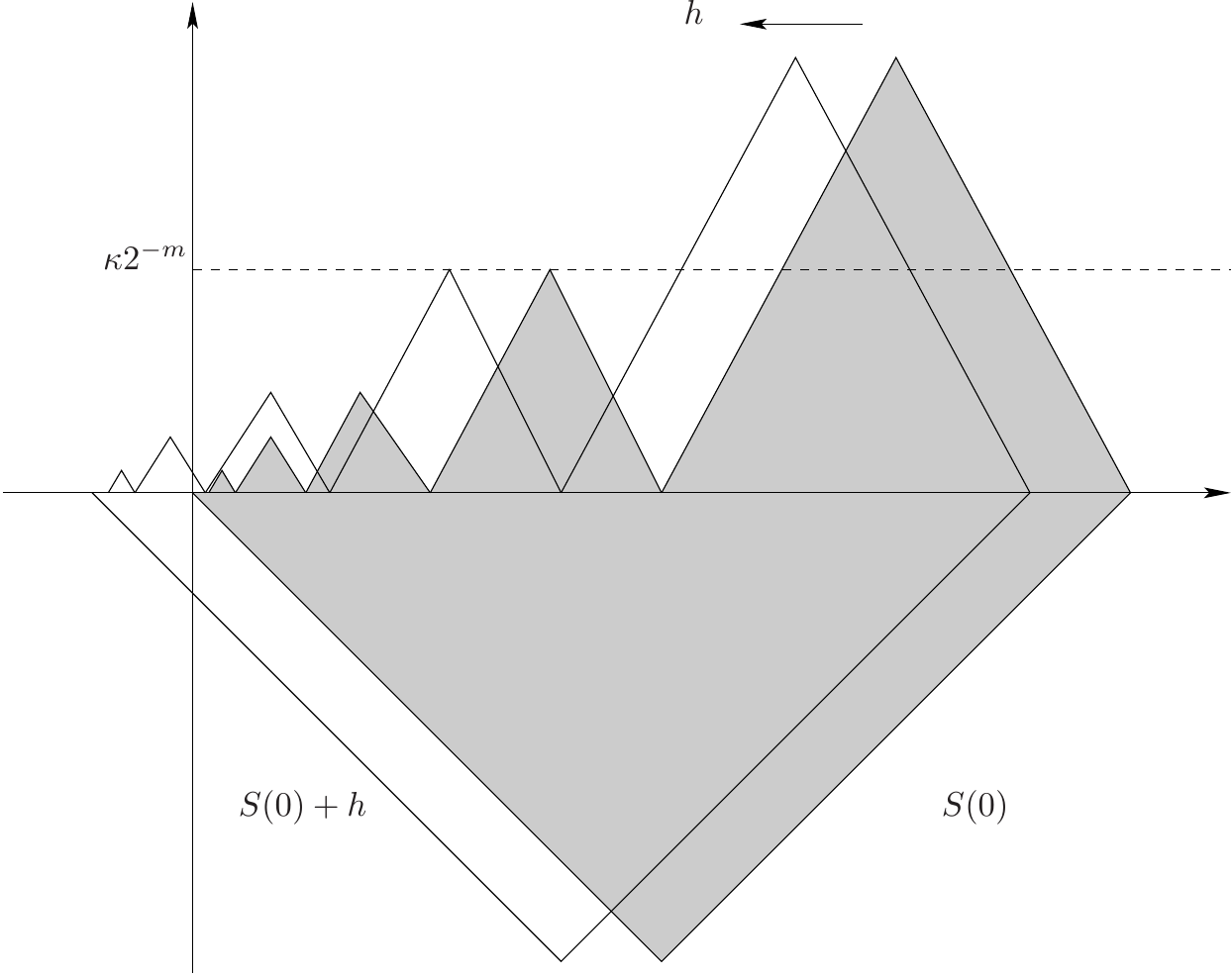} 
\caption{The choice of section line in Example 3.1.}
\end{figure}

\n
\smallskip

\n
We now show that for $\ve<1$, $h=\ve\,(-1,0)$ and $\eta=\eta\big((-1,0)\big)=(0,1)$, the sections of the symmetric difference $S_\Delta$ for all $y\in \Omega_\eta$ satisfy (\ref{eq:ex_H2_part1}).

\n Let $h=\ve(-1,0)$. Consider, as a preliminary step, the symmetric difference ${\cal T}_\Delta={\cal T}_m~\Delta~(h+{\cal T}_m)$ for any fixed $m\in\N$. It is not difficult to verify that for $y\in\Omega_\eta$
$$
\meas({\cal T}_\Delta)^y_\eta\leq 2\kappa\ve\,
$$
as illustrated by Figure 5. To pass to the sections of $S_\Delta$, it is sufficient to observe that  $(S_\Delta)_\eta$ is  the interval $(-\ve, 1)$ on the $x$--axis and that or any $y\in(S_\Delta)_\eta\subset\Omega_\eta$, there holds:
\begin{align*}
(S_\Delta)_\eta^y\subseteq \big({\cal T}_0~\Delta~(h+{\cal T}_0)\big)_\eta^y \cup 
\big({\cal T}_{\bar m}~\Delta~(h+{\cal T}_{\bar m})\big)_\eta^y,
\qquad\qquad &
\mbox{ when }~~
y\in(0,1)\,,\\
(S_\Delta)_\eta^y\subseteq \left\{
(a,b)\in\R^2~;~ a\in [-\ve,0]\,,\,b\in [-2\kappa a,2\kappa a]
\right\}_\eta^y,
\qquad\qquad &
\mbox{ when }~~y\in(-\ve,0] \,.
\end{align*}

\begin{figure}[h]
\includegraphics[width=.51\textwidth]{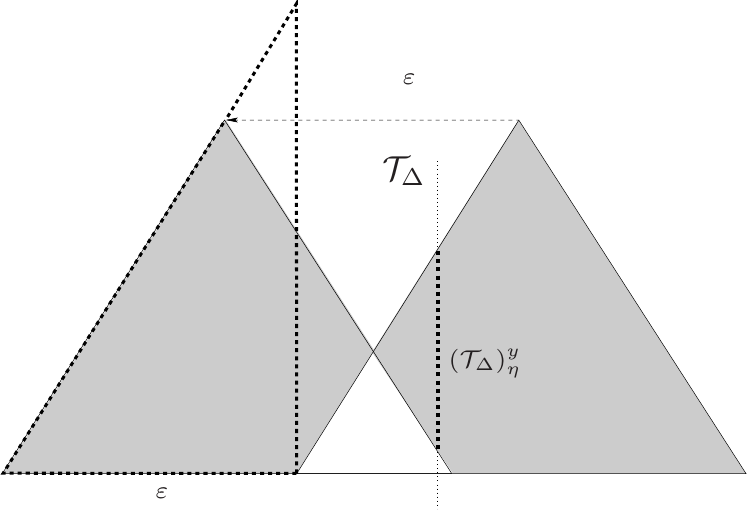}
~~~~~
\includegraphics[width=.40\textwidth]{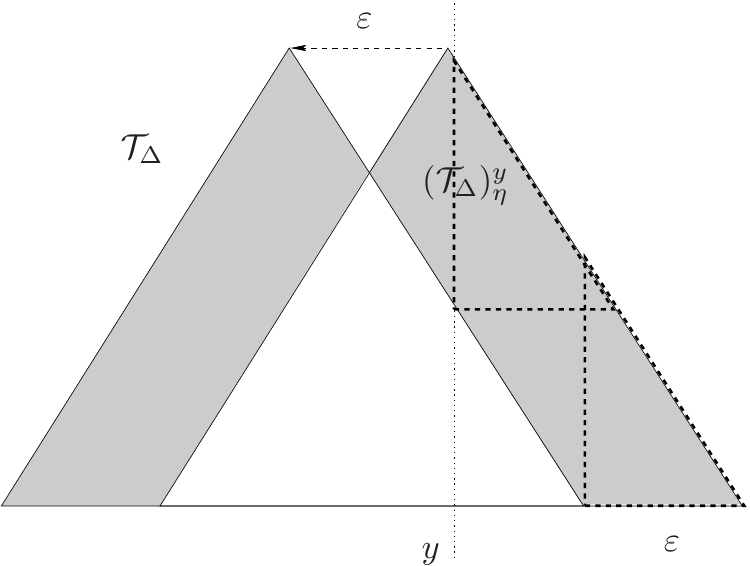}
\caption{
The measure of the sections $({\cal T}_\Delta)_\eta^y$ can be estimated by the height of the right triangle with
base $\ve$ and hypotenuse parallel to the side of ${\cal T}_m$. Left: $\ve\geq2^{-m-1}$. Right: $\ve<2^{-m-1}$.}
\label{fig:trianglesA}       
\end{figure}

\subsection{Swimmer's internal forces}\label{sec:force_model}

In this subsection we give  the precise description of the force term $ f$  in  (\ref{eq:nse}). Due to the nature of the swimming motion as self-propulsion, all the forces in  model (\ref{eq:nse})-(\ref{eq:ode}) are  \emph{internal} relative to the swimmer, i.e., their sum is equal to zero and their torque is constant. In turn, these forces, acting between swimmer's body parts, will  create a pressure upon the surrounding fluid and, thus, will act as external forces upon it.  We assume that all forces act through the immaterial links attached to the centers of mass $z_i(t)$'s of the sets $S(z_i(t))$, and then are uniformly transmitted to all points in their respective supports. 

\n
 Similar to~\cite{KhBook}, in this paper we consider two types of forces forming $ f$ in (\ref{eq:nse}): \emph{rotational forces} and \emph{elastic forces}, which we represent as
\begin{equation}\label{eq:forces}
f(t,x):= f_{rot}(t,x)+f_{el}(t,x)\,.
\end{equation}

\n The 2nd term  in (\ref{eq:forces}) describes the forces acting as elastic links between any two subsequent sets $ S (z_i)$'s to preserve the integrity  of swimmer's structure. They  act  according to the 3rd Newton's law and Hooke's law with variable (positive) rigidity coefficients 
$\kappa_1(t),\ldots,\kappa_{N-1}(t)$ when the distances between any two adjacent points $z_i(t)$ and $z_{i-1}(t)$, $i = 2,\ldots,N$, deviate from the respective given values $\ell_i>0, i = 1, \ldots, N-1$(see \cite{KhBook}): 
\begin{align}\label{eq:forces_el}
f_{el}(t,x)&:=\sum_{i=2}^{N} \bigg[\xi_{i-1}(t, x)\,\kappa_{i-1}(t)\,
{|z_i(t) - z_{i-1}(t)| -\ell_{i-1}\over |z_i(t) - z_{i-1}(t)|}\,
\big(z_i(t) - z_{i-1}(t)\big)\nonumber\\
&~~~~~~~+ \xi_{i}(t, x)\,\kappa_{i-1}(t)\,
{|z_{i}(t) -z_{i-1}(t)| - \ell_{i-1}\over |z_{i}(t) - z_{i-1}(t)|}\,
\big(z_{i-1}(t) - z_{i}(t)\big) \bigg].
\end{align}
Here
$\xi_i(t,x)$ denotes the characteristic function of   $S (z_i(t))$, i.e.,
$$
\xi_i(t,x) :=
\left\{
\begin{array}{ll}
1 & \mbox{if }x\in S (z_i(t))\\
0 & \mbox{otherwise }.
\end{array}
\right.\,
$$
\begin{rem}
In the above structure we can also assume that $\kappa_1(t),\ldots,\kappa_{N-1}(t)$ can be of any sign and replace  (\ref{eq:forces_el}) with more general (and simpler) description of pairs of co-linear forces between $ z_i$'s:
\begin{align}\label{eq:forces_el2}
f_{el}(t,x)&:=\sum_{i=2}^{N} \bigg[\xi_{i-1}(t, x)\,\kappa_{i-1}(t)\,
\big(z_i(t) - z_{i-1}(t)\big)\nonumber\\
&~~~~~~~+ \xi_{i}(t, x)\,\kappa_{i-1}(t)\,
\big(z_{i-1}(t) - z_{i}(t)\big) \bigg].
\end{align}
\end{rem}

\bigskip
\n
All the proofs in this paper are given for technically more complex case of elastic forces in (\ref{eq:forces_el}).

\bigskip

\n The 1st term in (\ref{eq:forces})  describes the forces that  allow each point $z_i(t), i = 2, \ldots, N-1$ to rotate any pair of the adjacent points $z_{i-1}(t)$ and $z_{i+1}(t)$ about it in either folding or unfolding fashion. In turn, by the 3rd Newton's law, these points will act back on $z_i(t)$ with the respective countering force. 
The description of rotational forces requires principally different approaches for the 2-$D$ and 3-$D$ cases.

\underline{In the 2-$D$ case} all the forces lie in the same plane at all times, and we can describe them by making use of the matrix
$$
A:=\left(\begin{array}{cc}0 & 1\\ -1 & 0\end{array}\right)\,
$$
as follows (see \cite{KhBook}): 
\begin{align}\label{eq:forces_rot2}
f_{rot,2d}(t,x)&:=\sum_{i=2}^{N-1} v_{i-1}(t)\bigg[\xi_{i-1}(t, x)\, A\,\big(z_{i-1}(t)-z_i(t)\big)\nonumber\\
&~~~~~~~-\xi_{i+1}(t, x)\,
{|z_{i-1}(t)-z_i(t)|^2\over |z_{i+1}(t)-z_i(t)|^2}\,
A\,\big(z_{i+1}(t)-z_i(t)\big)\bigg]\nonumber\\
&+ \sum_{i=2}^{N-1}
\xi_{i}(t, x)\,v_{i-1}(t)\bigg[A\,\big(z_{i}(t)-z_{i-1}(t)\big)\nonumber\\
&~~~~~~~-\,{|z_{i-1}(t)-z_i(t)|^2\over |z_{i+1}(t)-z_i(t)|^2}\,
A\,\big(z_{i}(t)-z_{i+1}(t)\big)\bigg]\,,
\end{align}
where functions $v_1(t),\ldots, v_{N-2}(t)$ characterize the strength and orientation (folding or unfolding) of  respective pairs of rotational forces at time $ t$.

In turn, \underline{in the 3-$D$ case}, to satisfy the 3rd Newton's law, we need to make sure that the respective rotational forces acting on $z_{i-1}(t)$ and $z_{i+1}(t)$ lie in the same plane spanned by the vectors $z_{i-1}(t) - z_i(t)$ and $z_{i+1}(t) -z_i (t)$. 
 In order to achieve the continuity of these forces in time, in this paper we choose to  reduce their magnitudes to zero, when the triplet  $\{z_{i-1}(t), z_i(t), z_{i+1}(t)\}$ approaches the aligned configuration (for other options  see \cite{Kh5}). Indeed, such  configuration admits infinitely many planes containing this triplet, which makes it  an intrinsic point of discontinuity for the procedure of the choice of the rotational  plane by means of the rotational forces whose magnitudes are strictly separated from zero. Respectively, we define the $3$d rotational forces as follows:
 \begin{align}\label{eq:forces_rot3}
f_{rot,3d}(t,x)&:=\sum_{i=2}^{N-1} v_{i-1}(t)\bigg[\xi_{i-1}(t, x)\, ~P_i[t]\big(z_{i-1}(t)-z_i(t)\big)\nonumber\\
&~~~~~~~-\xi_{i+1}(t, x)\,
{|z_{i-1}(t)-z_i(t)|^2\over |z_{i+1}(t)-z_i(t)|^2}\,
~Q_i[t]\big(z_{i+1}(t)-z_i(t)\big)\bigg]\nonumber\\
&+ \sum_{i=2}^{N-1}
\xi_{i}(t, x)\,v_{i-1}(t)\bigg[P_i[t]\big(z_{i}(t)-z_{i-1}(t)\big)\nonumber\\
&~~~~~~~-\,{|z_{i-1}(t)-z_i(t)|^2\over |z_{i+1}(t)-z_i(t)|^2}\,
~Q_i[t]\big(z_{i}(t)-z_{i+1}(t)\big)\bigg]\,,
\end{align}
where the scalar functions $v_1(t),\ldots, v_{N-2}(t)$ control the magnitudes of the rotational forces and determine whether they  act in folding or unfolding fashion, and
$$
{\bf x}\mapsto P_i[t]{\bf x}:=\big[(z_{i-1}(t)-z_i(t)) \times (z_{i+1}(t)-z_i(t))\big]\times {\bf x}\,,
$$
$$
{\bf x}\mapsto Q_i[t]{\bf x}:={\bf x}\times \big[(z_{i-1}(t)-z_i(t)) \times (z_{i+1}(t)-z_i(t))\big]\,.
$$
Note that $P_i[t]  {\bf x} =-Q_i[t] {\bf x}$ and  $|P_i[t]{\bf x}|=|Q_i[t]{\bf x}|\to 0$ for any  ${\bf x}$ when  points $z_{i-1}(t),z_i(t),z_{i+1}(t)$ converge to the  aligned configuration.

\bigskip
\n
\begin{rem}
The forcing term $f$  in (\ref{eq:nse}), \eqref{eq:forces} can (more precisely) be denoted as $f(t,x\,;z,\kappa,v)$.  However, we will use a shorter notation $f(t,x)$ as, typically, there is no ambiguity about the choices of $ z,\kappa$ and $v$. 
\end{rem}

\n
We can prove the following result.

\begin{lem}\label{lem:forces} Let $z=(z_1,\ldots,z_N) \in [\con([0,T]\,; \Omega )]^N \subset  [\con([0,T]\,;\R^d)]^N$, $\kappa=(\kappa_1,\ldots,\kappa_{N-1})\in [\Lsp^2(0,T)]^{N-1}$ and $v=(v_1,\ldots,v_{N-2})\in [\Lsp^2(0,T)]^{N-2}$ be fixed. Assume that for all $t\in[0,T]$ there holds
$$
|z_i(t)-z_j(t)| > 2r\,, \qquad\qquad i, j = 1, \ldots, N,~~ i \neq j\,,
$$
with $r>0$ as in {\bf (H1)}.
Then, the forcing term $f(t,x)$ defined in~\eqref{eq:forces} belongs to $\Lsp^2(0,T\,;[\Lsp^2(\Omega)]^d)$ and there hold the following estimates
\begin{align*}
\|f_{el}\|_{\Lsp^2(0,T\,;[\Lsp^2(\Omega)]^d)}~&\leq~2 \sqrt{\meas(\Omega)}\,
\|\sum_{i=2}^{N} |\kappa_{i-1}|\|_{\Lsp^2(0,T)}\\
&~~~~~~\times \max_{i=2,\ldots,N}\left\{
\|z_i-z_{i-1}\|_{\con([0,T]\,;\R^d)}+\ell_{i-1}
\right\}
\,,
\\
\|f_{rot,2d}\|_{\Lsp^2(0,T\,;[\Lsp^2(\Omega)]^2)}~&\leq~2\,\, \sqrt{\meas(\Omega)}\,
\|\sum_{i=2}^{N-1} |v_{i-1} |  \|_{\Lsp^2(0,T)}\\
&~~~~~~\times \max_{i=2,\ldots,N-1}\left\{
\|z_i-z_{i-1}\|_{\con([0,T]\,;\R^2)}+\,{\|z_i-z_{i-1}\|^2_{\con([0,T]\,;\R^2)}\over 2r}
\right\}
\,,
\\
\|f_{rot,3d}\|_{\Lsp^2(0,T\,;[\Lsp^2(\Omega)]^3)}~&\leq~4\,\, \sqrt{\meas(\Omega)}\,
 \|\sum_{i=2}^{N-1} |v_{i-1} | \|_{\Lsp^2(0,T)}\\
 &~~~~~~\times \max_{i=2,\ldots,N} \|z_i-z_{i-1}\|^3_{\con([0,T]\,;\R^3)}
\,.
\end{align*}
The above estimates yield that
\begin{equation}\label{z}
\| f \|_{\Lsp^2(0,T\,;[\Lsp^2(\Omega)]^d)} \; \leq \; \zeta,
\end{equation}
where positive constant $ \zeta > 0$ depends on $ T, \Omega$ and $ \Lsp^2(0,T)$-norms of parameters $ v_i$'s and $ \kappa_i$'s.
\end{lem}
\begin{rem}
It is not difficult to see that the above-defined forces in (\ref{eq:forces})-(\ref{eq:forces_rot3})  are internal relative to the swimmer at hand, i.e., their sum is zero and their angular momentums are constant (see, cf.~\cite[Chapter~12]{KhBook}).
\end{rem}


\section{Proofs of the main results}\label{sec:proof_main}
The proof of Theorem~\ref{thm:wellposed} will be based on Schauder's fixed point theorem, and will follow from  a series of lemmas. The general scheme of our proofs  is similar to that introduced  in \cite{Kh11}, \cite{Kh1}, \cite{KhBook} (Chapters 11-12), \cite{Kh5} for the case of the fluid governed by the linear nonstationary Stokes equations. In  this paper  we consider the full nonlinear Navier-Stokes equations, which will require  some principal technical modifications of this scheme. Our arguments below are nearly the same for  both 2-$D$ and 3-$D$ cases. However, in the part  dealing with the  Navier Stokes equations,  the  latter is traditionally  more challenging.

\bigskip
\n
\begin{rem}\label{CS} In our our proofs given below we assume that all the parameters in Theorems \ref{thm:wellposed} and \ref{thm:wellposed3D1} are  fixed, namely:  $ T>0$, the initial datum $u_0$ either in $H$ (if $d=2$) or in $V$ (if $d=3$), a vector  $z_0=(z_{1,0},\ldots,z_{N,0})\in\R^{dN}$ as in (\ref{body1})  and   functions $\kappa_1,\ldots,\kappa_{N-1}$, $v_1,\ldots,v_{N-2} \in \Lsp^2(0,T)$. {\sf We may also  omit the explicit mentioning of  dependence of some of the constants below from these  parameters}.
\end{rem}

\bigskip
\n
We will use the symbol $ T_1$ to denote the length of the time-interval of existence of solutions to Navier Stokes equations in (\ref{eq:nse}). For $ d = 2$ ,  due to  Theorem \ref{thm:existence2d} cited in Section 6, $T_1=+\infty $.

\bigskip
\n
For $ d = 3$, the size of $T_1 > 0$ depends on $(f, u_0)$ as cited in   Theorem~\ref{thm:existence3d} in Section 6 (Appendix A). 

More precisely, in this case we will  need to identify the value of  $T_1$ which, given $ u_0$,  will work  uniformly for any selection of  $f$ used in the proofs of Theorems  \ref{thm:wellposed}-\ref{thm:wellposed3D1}  below.
Namely, we  set:
\begin{equation}\label{T1}
T_1 \in (0,  \min \left\{T,  C^o \left[ ||  u_0 ||^2_{ H^1_0} + \zeta^2   \right]^{-2}\right\}) \;\;\; (d = 3),
\end{equation}
where $ \zeta$ is from (\ref{z}).

\n
In  subsection \ref{mainproof1} we introduce necessary notations and will state the aforementioned lemmas, whose proofs will be given in section \ref{sec:tech_proofs}. In  subsection \ref{mainproof2} we will prove Theorems \ref{thm:wellposed}-\ref{thm:wellposed3D1}.

\subsection{Formulation of auxiliary lemmas}\label{mainproof1}

\n
{\bf Set ${\cal B}_q$.} For  any  $T^*\in(0,T_1)$ and  $q>0$, set
\begin{align*}
{\cal B}_q&:=\left\{ \phi\in \Lsp^2(0,T^*\,;V)~;~\int_0^{T^*}\|\phi(t,\cdot)\|^2_{H^1_0}\,dt\leq q^2\right\}\,.\\
\end{align*}

\n
We will show below that, for  suitable choice of pair $ (T^*, q)$, the following three nonlinear operators will be well defined.

\bigskip

\n
{\bf Operator ${\bf T}$.} Set
\begin{equation}\label{eq:operator_ODE}
\begin{array}{c}
 {\bf T}\colon     \Lsp^2(0,T^*\,;V) \supset {\cal B}_q
~~\longrightarrow~[\con([0,T^*]\,;\Omega)]^N \subset [\con([0,T^*]\,;\R^d)]^N,\\\\
{\bf T}u := z = (z_1, \dots, z_N),
\end{array}
\end{equation}
where $z_1,\ldots,z_N$ are trajectories of system \eqref{eq:ode}  with the aforementioned fixed initial state $z_0$. We will show that for suitable pair $ (T^*, q)$ they will satisfy (\ref{body2}).

\bigskip
\n
{\bf Operator ${\bf F}$.} For any $ z \in {\bf T} ({\cal B}_q)$, set:
\begin{equation}\label{eq:operator_forces}
\begin{array}{c}
{\bf F}\colon[\con([0,T^*]\,;\R^d)]^N  \supset [\con([0,T^*]\,;\Omega)]^N \supset {\bf T} ({\cal B}_q) ~~ \longrightarrow ~~
\Lsp^2(0,T^*\,; [\Lsp^2(\Omega)]^d),\\\\
{\bf F}z={\bf F}(z_1, \dots, z_N):= f,
\end{array}
\end{equation}
where $f$ is the forcing term defined in~\eqref{eq:forces}, \eqref{eq:forces_rot2}-\eqref{eq:forces_rot3}.

\bigskip
\n
{\bf Operator ${\bf S}$.} For  any $ f \in {\bf F} {\bf T} ({\cal B}_q)$,  we set respectively:
\begin{equation}\label{eq:operator_NSE}
\begin{array}{c}
{\bf S}\colon \Lsp^2(0,T^*\,; [\Lsp^2(\Omega)]^d) \supset {\bf F} {\bf T} ({\cal B}_q)~~ \longrightarrow ~~
 \Lsp^2(0,T^*\,;V), \\\\
{\bf S}f:= u,
\end{array}
\end{equation}
where $u$ is the solution to uncoupled Navier--Stokes equations~\eqref{eq:nse}.

\bigskip
\n
To prove Theorem 1, we intend to show that the mapping 
$$
 {\bf S}{\bf F}{\bf T}\colon{\cal B}_q
~~\longrightarrow~~  \Lsp^2(0,T^*\,;V)
$$
has a unique fixed point in ${\cal B}_q$ for some pair $(T^*,q)$. To this end, we will show that all the above operators  are continuous, while ${\bf T}$ is also compact. Then the desirable result will follow by Schauder's fixed point theorem.

\begin{lem}\label{lem:ode} Fix $T>0$, 
$u\in L^2(0,T\,; [H^1_0(\Omega)]^d)$, $ d = 2,3$, and assume (\ref{body1}).
Then there exists $T^*\in(0,T]$ 
such that system~\eqref{eq:ode} admits a unique solution $z$ in 
$[\con([0,T^*]\,;\Omega)]^N \subset  [\con([0,T^*]\,;\R^d)]^N$ and (\ref{body2}) holds.
Finally, given  a constant $k>0$, there exists $T^*_k\in(0,T]$ such that (\ref{body2})  holds on $[0,T^*_k]$  (in place of $ [0, T^*]$) for any choice of  $u$ such that $\|u\|_{L^2(0,T\,; [H^1_0(\Omega)]^d)}\leq k$.
\end{lem}

\begin{lem}\label{lem:T} Under the assumptions of Lemma \ref{lem:ode}, for any fixed $q>0$, let $T^*$ be such that
\begin{equation}\label{eq:time_traj_T}
0< T^* < \min \left\{
T^*_q,  \, \left({\meas(S(0))\over \,c_\Omega\,{c'_\Omega}\,{\cal K}_S\, q }\right)^2
\right\} ,
\end{equation}
where  $T^*_q$ is defined in  Lemma~\ref{lem:ode} (see (\ref{Tq})), $c_\Omega$ is the constant in~\eqref{eq:sobolev_sec}, $c'_\Omega$ is the constant in~\eqref{eq:fubini_sec}  and ${\cal K}_S$ is the constant in {\bf (H2)}.
Then, operator ${\bf T}$ defined in~\eqref{eq:operator_ODE} is well defined, bounded, continuous and compact. In particular,
there exist constants $L_{\bf T},L'_{\bf T}>0$ such that for all $u,v\in{\cal B}_q$:
\begin{equation}\label{eq:T_bdd}
\|{\bf T}u\|_{[\con([0,T^*]\,;\R^d)]^N}
\leq L_{\bf T}\Big(|z_0| +
\,\sqrt{T^*}\,\|u\|_{ \Lsp^2(0,T^*\,;V)}\Big)\,,
\end{equation}
and
\begin{equation}\label{eq:T_cont}
\|{\bf T}u-{\bf T}v\|_{[\con([0,T^*]\,;\R^d)]^N}
\leq L_{\bf T}'
\,{\sqrt{T^*}\over 1-\,{c_\Omega\,c'_\Omega\,{\cal K}_S\,q\,\sqrt{T^*}\over \meas(S(0))}}
\,\|u-v\|_{ \Lsp^2(0,T^*\,;V)}\,.
\end{equation}
\end{lem}

\begin{lem}\label{lem:F} For any fixed  $ q >0 $ and $T^*$ as in Lemma \ref{lem:T}, 
operator ${\bf F}$  in~\eqref{eq:operator_forces} is well defined, bounded and continuous. In particular,
there exist positive constants $L_{\bf F}, L'_{\bf F} $ such that for 
$
z,w\in{\bf T} ({\cal B}_q)$ we have:
\begin{equation}\label{eq:F_bdd}
\|{\bf F}z\|_{\Lsp^2(0,T^*\,;[\Lsp^2(\Omega)]^d)}
\leq L_{\bf F}\,\Gamma_{\kappa,v}(T^*) \leq \zeta,
\end{equation}
where the last inequality is due to (\ref{z}), and
\begin{equation}\label{eq:F_cont}
\|{\bf F}z-{\bf F}w\|_{\Lsp^2(0,T^*\,;[\Lsp^2(\Omega)]^d)}
\leq L'_{\bf F}\,\Gamma_{\kappa,v}(T^*)
~\|z-w\|_{[\con([0,T^*]\,;\R^d)]^N}\,,
\end{equation}
$$
\Gamma_{\kappa,v}(T^*) :=\max\left\{\|\kappa_1\|_{\Lsp^2(0,T^*)},\ldots,\|\kappa_{N-1}\|_{\Lsp^2(0,T^*)},\|v_1\|_{\Lsp^2(0,T^*)},\ldots,\|v_{N-2}\|_{\Lsp^2(0,T^*)}\right\}\,.
$$
\end{lem}
\begin{rem}\label{4.1} The above estimates (\ref{eq:F_bdd})-(\ref{eq:F_cont})  can be refined to include explicitly 
 the $ [\con([0,T^*]\,;\R^d)]^N$-norms of trajectories $ z,w\in{\bf T} ({\cal B}_q)$. However, this is not principal, because   $ {\bf T} ({\cal B}_q) \subset [\con([0,T^*]\,;\Omega)]^N$.
\end{rem}

\begin{lem}\label{lem:S}  For any fixed $T^*\in(0,T_1)$, operator ${\bf S}$  defined in~\eqref{eq:operator_NSE} is well defined, bounded and continuous. In particular,
there exist constants $L_{\bf S},L_{\bf S}'$ and $C_{\bf S} $ such that (see Remark \ref{CS})
\begin{equation}\label{eq:S_bdd}
\|{\bf S} f\|_{ \Lsp^2(0,T^*\,;V)} 
\leq L_{\bf S}
\left(\|u_0\|_{\Lsp^2}\,+\,\|f\|_{\Lsp^2(0,T^*\,;[\Lsp^2(\Omega)]^d)}\right)\,
\end{equation}
and for $ d = 2$:
$$
\|{\bf S}f-{\bf S}g\|_{ \Lsp^2(0,T^*\,;V)}
$$
\begin{equation}\label{eq:S_cont2d}
\leq L_{\bf S}' \, e^{C_{\bf S}\, \max\{\|f\|^2_{ \Lsp^2(0,T^*\,;[\Lsp^2(\Omega)]^2)},\|g\|^2_{ \Lsp^2(0,T^*\,;[\Lsp^2(\Omega)]^2)}\}}
\,\|f-g\|_{\Lsp^2(0,T^*\,;[\Lsp^2(\Omega)]^2)}\,,
\end{equation}
while for $ d = 3$:
$$
\|{\bf S}f-{\bf S}g\|_{ \Lsp^2(0,T^*\,;V)}
$$
\begin{equation}\label{eq:S_cont3d}
\leq L_{\bf S}' \, e^{C_{\bf S}  \, \max\{\|f\|^8_{ \Lsp^2(0,T^*\,;[\Lsp^2(\Omega)]^3)},\|g\|^8_{ \Lsp^2(0,T^*\,;[\Lsp^2(\Omega)]^3)}\}}
\,\|f-g\|_{\Lsp^2(0,T^*\,;[\Lsp^2(\Omega)]^3)}\,.
\end{equation}
\end{lem}


\subsection{Proof of Theorems Theorems~\ref{thm:wellposed}-\ref{thm:wellposed3D1}}\label{mainproof2}

\n
{\bf Existence.} Consider any $ q>0$. Then, 
Lemmas \ref{lem:ode}- \ref{lem:S} imply that operator ${\bf SFT}$,
\begin{equation}\label{oppr}
{\bf SFT}\colon {\cal B}_q
~~ \longrightarrow ~~
\Lsp^2(0,T^*\,;V)
\end{equation}
is well defined, continuous and compact for any choice of $T^*$ satisfying~\eqref{eq:time_traj_T}.

\n
We now claim that, by suitable choices of $q$ and further refinement (i.e., possible reducing) of the value of $T^*$, we can prove that 
${\bf SFT}({\cal B}_q)\subseteq {\cal B}_q$.

\n
Let us choose
$q\geq 2L_{\bf S}\|u_0\|_{\Lsp^2}$, where $ L_{\bf S}$ is from   (\ref{eq:S_bdd}).  Then, taking into account ~\eqref{eq:time_traj_T}, choose
$$
0<T^*<\min\left\{1\,,\, T^*_q\,,\,T^*_{\kappa,v}\,,\,\left({\meas(S(0))\over c_\Omega\,c'_\Omega\,{\cal K}_S\,q}\right)^2\right\}\,,
$$
where  $T^*_{\kappa,v}>0$  is small enough to ensure that  
\begin{equation}\label{eq:time_small_enough}
L_{\bf S}  L_{\bf F} \Gamma_{\kappa,v}(T^*_{\kappa,v}) \leq\, {q\over 2}\,, 
\end{equation}
see    Lemmas  \ref{lem:F}-\ref{lem:S} for parameters $L_{\bf F},  \Gamma_{\kappa,v}(T^*_{\kappa,v})$.
Then,  making use of  (\ref{eq:S_bdd}) and (\ref{eq:F_bdd}),  we obtain:
\begin{align*}
\|{\bf SFT} u\|_{ \Lsp^2(0,T^*\,;V)} & \leq 
L_{\bf S}\,\|u_0\|_{\Lsp^2}+ L_{\bf S}\,\|{\bf FT}u\|_{\Lsp^2(0,T^*\,;[\Lsp^2 (\Omega)]^d)}\\
&\leq \,{q\over 2}\,+L_{\bf S}\,L_{\bf F}\, \Gamma_{\kappa,v}(T^*) \, \leq q\,.
\end{align*}

\n Applying Schauder's fixed point theorem completes the proof of existence in Theorems~\ref{thm:wellposed} and \ref{thm:wellposed3D1}.~~$\diamond$ 

\medskip

\n 
{\bf Proof of uniqueness.} We argue by contradiction. 
 Let $(u^{(1)},z^{(1)})$ and $(u^{(2)}, z^{(2)})$ be two distinct solutions to~\eqref{eq:nse}, \eqref{eq:ode} for the same initial condition and forcing term, both defined on the same time-interval $[0,\hat{T}]$. Due to linearity of (\ref{eq:ode}) and Lemma \ref{lem:ode}, this means that $u^{(1)} $ and $ u^{(2)}$ must be  distinct  on $[0,\hat{T}]$. Hence, since  solutions  to ~\eqref{eq:nse} are continuous in time with values in $H$, there exists some  interval, say,  $(T,T+\tau), T \geq0, \tau > 0$, such that   these two solutions are distance at every point of it (in the above sense of continuity) but  $u^{(1)}(T, \cdot) = u^{(2)}(T, \cdot)= u(T)$. Setting (similar to the proof of existence):
$$
\hat q:=\max\left\{\|u^{(1)}\|_{ \Lsp^2(T,T^{(1)}\,;V)}\,, \|u^{(2)}\|_{ \Lsp^2(T,T^{(1)}\,;V)}\,, 2L_{\bf S}\|u^{(1)}(T)\|_{\Lsp^2}\right\},
$$
we choose $\tau^*\in(0,T^{(1)}-T]$ so that 
$$
0<\tau^*<\min\left\{1\,,\, T^*_{\hat q}\,,\,\hat T^*_{\kappa,v}\,,\,
\left({\meas(S(0))\over 2\,c_\Omega\,c'_\Omega\,{\cal K}_S\,\hat q}\right)^2\right\},
$$
with $\hat T^*_{\kappa,v}$ as in~\eqref{eq:time_small_enough} with $\hat q/2$ in place of $q/2$.
In particular, by defining the operators ${\bf S},{\bf F},{\bf T}$ as in~\eqref{eq:operator_ODE}--\eqref{eq:operator_NSE} with the time interval $[0,T^*]$ replaced by $[T,T+\tau^*]$ and with ${\cal B}_q$ replaced by ${\cal B}_{\hat q}\subset \Lsp^2(T,T+\tau^*\,;V)$, we can repeat the previous calculations yielding  ${\bf SFT}{\cal B}_{\hat q}\subseteq{\cal B}_{\hat q}$ and that Lemmas~\ref{lem:T}--\ref{lem:S} apply.

\n Since both $u^{(1)}$ and $u^{(2)}$ are fixed points of ${\bf SFT}$ on ${\cal B}_{\hat q}$, we have from 
(\ref{eq:S_cont2d}) in the case when $ d = 2$ (and we use (\ref{eq:S_cont3d}) for $ d = 3$):
$$
 \|u^{(1)}-u^{(2)}\|_{ \Lsp^2(T, T+\tau^*\,;V)}  = \|{\bf SFT} u^{(1)}-{\bf SFT} u^{(2)}\|_{ \Lsp^2(T, T+\tau^*\,;V)} 
$$
$$
\leq  L'_{\bf S} \,e^{C_{\bf S} \max_{j=1,2}  \|{\bf FT}u^{(j)}\|^2_{\Lsp^2(T, T+\tau^*\,; [\Lsp^2(\Omega)]^d)}} \,
\|{\bf FT}u^{(1)}-{\bf FT} u^{(2)}\|_{\Lsp^2(T, T+\tau^*\,;[\Lsp^2 (\Omega)]^d)}
$$
In turn, from (\ref{eq:F_cont}) it follows that$$
\|{\bf FT}u^{(1)}-{\bf FT} u^{(2)}\|_{\Lsp^2(T, T+\tau^*\,;[\Lsp^2]^d)} 
$$
$$
\leq 
L'_{\bf F}\,\Gamma_{\kappa,v}(\tau^*)\,
\,\|{\bf T}u^{(1)}-{\bf T}u^{(2)}\|_{[\con([T, T+\tau^*]\,;\R^d)]^N}.
$$

\n
Finally,  (\ref{eq:T_cont}) yields:
$$
\|{\bf T}u^{(1)}-{\bf T}u^{(2)}\|_{[\con([T, T+\tau^*]\,;\R^d)]^N} \leq  L_{\bf T}'
\,{\sqrt{\tau^*}\over 1-\,{c_\Omega\,c'_\Omega\,{\cal K}_S\, \hat q\,\sqrt{\tau^*}\over \meas(S(0))}}
\,\|u^{(1)}-u^{(2)}\|_{ \Lsp^2(T, T+\tau^*\,;V)}.
$$
Combining the above chain of estimates  allows us to conclude that
$$
\|u^{(1)}-u^{(2)}\|_{ \Lsp^2(T, T+\tau^*\,;V)}  \leq \Upsilon (\tau^*) 
\,\|u^{(1)}-u^{(2)}\|_{ \Lsp^2(T, T+\tau^*\,;V)}\,,
$$
where  $\Upsilon (\tau^*) \rightarrow 0$ as $ \tau^*$ tends to zero. 
Thus, $u^{(1)}$ and $u^{(2)}$  must  coincide on  $ [T,T+\tau_*]$ for some $  \tau_* \in (0, \tau^*)$. Contradiction.~~$\diamond$

\bigskip
\n


\section{Proofs of Lemmas \ref{lem:ode}- \ref{lem:S}}\label{sec:tech_proofs}

Below we use the symbol $ | a |, a \in \R^d$ to denote $ || a ||_{\R^d}$.

\subsection{Proof of Lemma~\ref{lem:ode}}

Fix $T>0$ and $u\in L^2(0,T\,; [H^1_0(\Omega)]^d)$ and assume (\ref{body1}). 

\n
{\bf Step 1.} Let $h_0$ and ${\cal K}_S$ as in {\bf (H2)}, and choose $T_0$ such that
\begin{equation}\label{eq:T0}
0<T_0<\min\left\{
T,\,{\meas(S(0))\,h_0^2\over 4\,\|u\|^2_{\Lsp^2(0,T\,;[\Lsp^2 (\Omega)]^2)}}\,,\,\left({\meas(S(0))\over 2\,c_\Omega\,{c'_\Omega}\,{\cal K}_S\,\|u\|_{\Lsp^2(0,T\,;[H^1_0 (\Omega)]^2)}}\right)^2
\right\},
\end{equation}
where $c_\Omega$ is the constant in~\eqref{eq:sobolev_sec} and $c'_\Omega$ is the constant in~\eqref{eq:fubini_sec}. Denoting
$$
B_{h_0/2}:=\left\{ 
\zeta\in\con([0,T_0]\,;\R^d)~;~\|\zeta\|_{\con([0,T_0]\,;\R^d)}\leq\,{h_0\over 2}
\right\},
$$
we define, for each index $i=1,\ldots,N$, a function $D_i\colon z_{i,0}+B_{h_0/2}\to\con([0,T_0]\,;\R^d)$ by setting for all $w\in z_{i,0}+B_{h_0/2}$
$$
D_i(w)(t):=z_{i,0}+\,{1\over\meas(S(0))}\,\int_0^t\int_{S(w(\tau))}\!\!u(\tau,x)\,dx\,d\tau\,, \;\; t \in [0, T_0].
$$

\medskip

\n{\bf Step 2.}  We claim that each $D_i$ takes values in $z_{i,0}+B_{h_0/2}$ and is a contraction mapping.

\n
{\it Proof of 
claim.} For all $t\in [0,T_0]$, we have, :
\begin{align*}
|D_i(w)(t)-z_{i,0}|& \leq \,{1\over\meas(S(0))}\,| \int_0^t\int_{S(w(\tau))} \!\! u(\tau,x) \,dx\,d\tau |\\
&\leq \, \,{\sqrt{T_0}\over\sqrt{\meas(S(0))}}\,\|u\|_{\Lsp^2(0,T\,;[\Lsp^2 (\Omega)]^d)}<\,{h_0\over 2}\,,
\end{align*}
where we have used {\bf (H1)}, Schwartz  inequality and~\eqref{eq:T0}. This proves that  $D_i\colon z_{i,0}+B_{h_0/2}\to z_{i,0}+B_{h_0/2}$.

\n Moreover, if we take $w_1,w_2\in z_{i,0}+B_{h_0/2}$, then  for all $t\in [0,T_0]$:
\begin{align*}
| D_i(w_1)(t)-D_i(w_2)(t) |& \leq \,{1\over\meas(S(0))}\, \left| \int_0^t  \int_{S(w_1(\tau))} \!\! u(\tau,x)\,dx - \int_{S(w_2(\tau))} \!\! u(\tau,x) \,dx\,\,d\tau \right| \\
& = \,{1\over\meas(S(0))}\,  \left|\int_0^t \int_\Omega u(\tau,x)\,\big( \xi_{1,\tau}(x)-\xi_{2,\tau}(x)\big)\,dx\,\,d\tau \right|,
\end{align*}
where we have denoted with $\xi_{i,\tau}(x)$ the characteristic functions of $S(w_i(\tau))$. 

\n For all $\tau\in[0,t]$ such that  $w_1(\tau)\neq w_2(\tau)$, by setting $h(\tau):=w_1(\tau)-w_2(\tau)$ and denoting with $\nu(\tau):=\eta (h(\tau))$ the unit vector given by {\bf (H2)} for $h(\tau)\in B_{h_0}(0)\setminus \{0\}$, we can consider the sections of $u(\tau,\cdot)(\xi_{1,\tau}-\xi_{2,\tau})$ corresponding to $y\in\Omega_{\nu(\tau)}$ ($\Omega_{\nu(\tau)}$ is the projection of $\Omega$  on the hyperplane orthogonal to $\eta$). For these sections, there holds, see notations (\ref{eq:section_set}),
$$
\int_\Omega u(\tau,x)\,\big( \xi_{1,\tau}(x)-\xi_{2,\tau}(x)\big)\,dx=
\int_{\Omega_{\nu(\tau)}}\Bigg(\int_{\Omega_{\nu(\tau)}^y}
u_{\nu(\tau)}^y(\tau,\alpha)~\left(\xi_{1,\tau}-\xi_{2,\tau}\right)_{\nu(\tau)}^y(\alpha)\,d\alpha\Bigg)\,d{\cal H}^{d-1}(y),
$$
where ${\cal H}^{n-1}$ stands for the $(n-1)$-dimensional Hausdorff measure over Borel sets of $\R^n$. 
Then, due to ~\eqref{eq:sobolev_sec}, for each component $ u_j $ of $ u = (u_1, \dots, u_d)$:
\begin{align*}
\bigg|\int_\Omega u_j(\tau,x)\,&\big( \xi_{1,\tau}(x)-\xi_{2,\tau}(x)\big)\,dx\,\bigg|\\
&\leq\int_{\Omega_{\nu(\tau)}}\!\!\!\!\!\|(u_j)_{\nu(\tau)}^y(\tau,\cdot) \|_{\Lsp^\infty (\Omega_{\nu(\tau)}^y)}
~\Bigg(\int_{\Omega_{\nu(\tau)}^y}
\left|\left(\xi_{1,\tau}-\xi_{2,\tau}\right)_{\nu(\tau)}^y(\alpha)\right|\,d\alpha\Bigg)\,d{\cal H}^{d-1}(y)\\
&\leq \,c_\Omega\,\int_{\Omega_{\nu(\tau)}}\!\!\!\!\!\|(u_j)_{\nu(\tau)}^y(\tau,\cdot) \|_{H^1_0 (\Omega_{\nu(\tau)}^y)}
~\Bigg(\int_{\Omega_{\nu(\tau)}^y}
\left|\left(\xi_{1,\tau}-\xi_{2,\tau}\right)_{\nu(\tau)}^y(\alpha)\right|\,d\alpha\Bigg)\,d{\cal H}^{d-1}(y).
\end{align*}

\n Now notice that, $( \xi_{1,\tau} -\xi_{2,\tau})$ coincides with the characteristic function of the symmetric difference $S(w_1(\tau))\,\Delta\, S(w_2(\tau))$, and that we have
$$
S(w_1(\tau))\,\Delta\, S(w_2(\tau))~=~
\big(w_1(\tau)+S(0)\big)\,\Delta\, \big(w_2(\tau)+S(0)\big)~=~
\big(h(\tau)+S(0)\big)\,\Delta\, S(0),
$$
thanks to {\bf (H1)}. Thus, the above integral over $\Omega_{\nu(\tau)}^y$ is actually the measure of the section along $\eta(h(\tau))$ of this symmetric difference and, owing to {\bf (H2)}, we obtain:
\begin{equation}\label{eq:estimate_symm_diff}
\int_{\Omega_{\nu(\tau)}^y}
\left|\left(\xi_{1,\tau}-\xi_{2,\tau}\right)_{\nu(\tau)}^y(\alpha)\right|\,d\alpha
\leq {\cal K}_S\,|w_1(\tau)-w_2(\tau)|\leq {\cal K}_S\,\|w_1-w_2\|_{\con([0,T_0]\,;\R^d)}\,.
\end{equation}
Observing that, for times $\tau\in[0,t]$ such that $w_1(\tau)=w_2(\tau)$, the symmetric difference is empty and $\xi_{1,\tau}-\xi_{2,\tau}\equiv 0$ and, making use of~\eqref{eq:fubini_sec}, we can  deduce next that

$$
|D_i(w_1)(t)-D_i(w_2)(t)|^2 
 \leq 
 \,{1\over(\meas(S(0)))^2}\,  \left|\int_0^t \int_\Omega u(\tau,x)\,\big( \xi_{1,\tau}(x)-\xi_{2,\tau}(x)\big)\,dx\,\,d\tau \right|^2
$$
 $$ =  \,{1\over(\meas(S(0)))^2}\,  \sum_{j=1}^d \left|\int_0^t \int_\Omega u_j(\tau,x)\,\big( \xi_{1,\tau}(x)-\xi_{2,\tau}(x)\big)\,dx\,\,d\tau \right|^2 
 $$
$$ \leq   \,{1\over(\meas(S(0)))^2}\,  \sum_{j=1}^d  \{\int_0^t \int_\Omega 
 \,c_\Omega\,\int_{\Omega_{\nu(\tau)}}\!\!\!\!\!\|(u_j)_{\nu(\tau)}^y(\tau,\cdot) \|_{H^1_0 (\Omega_{\nu(\tau)}^y)}
 $$
 $$
 \times 
~\Bigg(\int_{\Omega_{\nu(\tau)}^y}
\left|\left(\xi_{1,\tau}-\xi_{2,\tau}\right)_{\nu(\tau)}^y(\alpha)\right|\,d\alpha\Bigg)\,d{\cal H}^{d-1}(y)\,dx\,\,d\tau \}^2 
$$
$$
\leq \,\{{c_\Omega\,{\cal K}_S\,\|w_1-w_2\|_{\con([0,T_0]\,;\R^d)}\over\meas(S(0))}\}^2 \, \sum_{j=1}^d  
\left[ \int_0^t\int_{\Omega_{\nu(\tau)}}\!\!\!\!\!\|(u_j)_{\nu(\tau)}^y(\tau,\cdot) \|_{H^1_0 (\Omega_{\nu(\tau)}^y)}\,d{\cal H}^{d-1}(y)\,d\tau\right]^2 
$$
$$
\leq \,\{{c_\Omega\,{\cal K}_S\,\|w_1-w_2\|_{\con([0,T_0]\,;\R^d)}\over\meas(S(0))}\}^2 \, \sum_{j=1}^d   \left[ \int_0^{T_0}  c_\Omega^\prime  \|u_j (\tau,\cdot) \|_{H^1_0 (\Omega)} d\tau\right]^2 
$$
$$
\leq \,\{{c_\Omega\,c'_\Omega\, {\cal K}_S\,\|w_1-w_2\|_{\con([0,T_0]\,;\R^d)}\over\meas(S(0))}\}^2 \,  T_0  \sum_{j=1}^d  \left[ \int_0^{T_0} \| u_j  (\tau,\cdot) \|^2_{H^1_0 (\Omega)}d\tau\right] 
$$
\begin{equation}\label{eq:estimate_with_slices}
\leq \,\{{c_\Omega\,c'_\Omega\, {\cal K}_S\,\|w_1-w_2\|_{\con([0,T_0]\,;\R^d)}\over\meas(S(0))}\}^2 \,  T_0 ,
\|u\|_{\Lsp^2(0,T_0\,;[H^1_0 (\Omega)]^d)}\,.
\end{equation}
Therefore, taking the supremum over $t\in [0,T_0]$ and recalling~\eqref{eq:T0}, we conclude that
$$
\|D_i(w_1)-D_i(w_2)\|_{\con([0,T_0]\,;\R^d)}<\,{1\over 2}\,\|w_1-w_2\|_{\con([0,T_0]\,;\R^d)}\,,
$$
and the proof of the above claim is completed.

\medskip

\n{\bf Step 3.} By applying the contraction mapping theorem to each $D_i$, we obtain uniquely determined fixed points $z^*_1,\ldots,z^*_N$ in $z_{i,0}+B_{h_0/2}\subset\con([0,T_0]\,;\R^d)$. Since  the choice of $T_0$ in~\eqref{eq:T0} is time invariant, we can repeat the same argument for intervals $[n\,T_0,(n+1)\,T_0]$, for $n\geq 1$, until we reach $T$. Thus, we obtain curves $z^*_1,\ldots,z^*_N$ in $\con([0,T]\,;\R^d)$ which satisfy the integral representations
$$
z^*_i(t)=z_{i,0}+\,{1\over\meas(S(0))}\,\int_0^t\int_{S(z^*_i(\tau))}\!\!u(\tau,x)\,dx\,d\tau\,,
\qquad\qquad t\in [0,T]\,,~~ i=1,\ldots,N\,,
$$
and, hence, they are solutions to~\eqref{eq:ode}.
By continuity in time of $z^*_i(t)$'s, it follows from (\ref{body1}) that we can choose $T^*\in (0,T]$ small enough to have
$$
\overline{S}(z^*_i(t)) \subset \Omega\,,
\qquad\qquad
|z^*_i(t)-z^*_j(t)| > 2r\,
$$
for all required indeces  and for all $t\in[0,T^*]$, so that~\eqref{body2} is satisfied. 

\medskip

\n{\bf Step 4.} 
Recalling that $\|u\|_{\Lsp^2(0,T_0\,;[\Lsp^2(\Omega)]^d)}\leq C_\Omega \|u\|_{\Lsp^2(0,T_0\,;[H^1_0]^d)}$, due to Poincar\'e inequality, it is immediate to see that the last statement in Lemma \ref{lem:ode} follows by replacing~\eqref{eq:T0} with 
\begin{equation}\label{Tq}
0<T_k<\min\left\{
T,\,{\meas(S(0))\,h_0^2\over 4\,C_\Omega\, k^2}\,,\left({\meas(S(0))\over 2\,c_\Omega\,{c'_\Omega}\,{\cal K}_S\, k}\right)^2
\right\} \,.
\end{equation}
This concludes the proof.~~$\diamond$

\subsection{Proof of Lemma~\ref{lem:forces}}

Let us start from the term $f_{el}$. It is immediate to see that, for any fixed $(t,x)\in (0,T)\times\Omega$, there holds
\begin{align*}
|f_{el}(t,x)|~&\leq~ \sum_{i=2}^{N} |\kappa_{i-1}(t)|\,\Big||z_i(t) - z_{i-1}(t)| -\ell_{i-1}\Big|\, \big(\xi_{i-1}(t, x)+\xi_i(t, x)\big)\\
&\leq~ 2\,\max_{i=2,\ldots,N}\left\{
\|z_i-z_{i-1}\|_{\con([0,T]\,;\R^d)}+\ell_{i-1}
\right\}\,\left(\sum_{i=2}^{N} |\kappa_{i-1}(t)|\right)\,.
\end{align*}
Thus, the first inequality in Lemma~\ref{lem:forces} follows.

\n
Similarly, for $f_{rot,2d}$ we have
\begin{align*}
|f_{rot,2d}(t,x)|~&\leq~ \sum_{i=2}^{N-1} |v_{i-1}(t)|\,|z_i(t) - z_{i-1}(t)|\, \big(\xi_{i-1}(t, x)+\xi_i(t, x)\big)\\
&~~~~~~~~+ \sum_{i=2}^{N-1} |v_{i-1}(t)| \,{|z_{i-1}(t)-z_i(t)|^2\over |z_{i+1}(t)-z_i(t)|}\,\big(\xi_i(t, x)+\xi_{i+1}(t, x)\big)\\
&\leq~ 2\, \max_{i=2,\ldots,N-1}\left\{
\|z_i-z_{i-1}\|_{\con([0,T]\,;\R^d)}+\,{\|z_i-z_{i-1}\|^2_{\con([0,T]\,;\R^d)}\over 2r}
\right\}\,\left(\sum_{i=2}^{N-1} |v_{i-1}(t)|\right)\,,
\end{align*}
because ${\bf x}\mapsto A{\bf x}$ is an isometry. Hence, for every $t\in(0,T)$ we deduce
\begin{align*}
\left(\int_\Omega |f_{rot,2d}(t,x)|^2\,dx\right)^{1/2}\!\!&\leq 2\,\left(\sum_{i=2}^{N-1} |v_{i-1}(t)|\right)~\sqrt{\meas(\Omega)}\\
&~~~~\times~
 \max_{i=2,\ldots,N-1}\left\{
\|z_i-z_{i-1}\|_{\con([0,T]\,;\R^d)}+\,{\|z_i-z_{i-1}\|^2_{\con([0,T]\,;\R^d)}\over 2r}
\right\}\,,
\end{align*}
so that the second inequality in Lemma~\ref{lem:forces} follows as well. 

\n
It remains to prove the third inequality about $f_{rot,3d}$. It is immediate to observe that for $i=2,\ldots,N-1$ there hold:
$$
\Big|P_i[t](z_{i-1}(t)-z_i(t))\Big|=\Xi_i(t)\cdot|z_{i-1}(t)-z_i(t)| \,,
\qquad
\Big|Q_i[t](z_{i+1}(t)-z_i(t))\Big|=\Xi_i(t)\cdot|z_{i+1}(t)-z_i(t)|\,,
$$
where $\Xi_i(t):=\big|(z_{i-1}(t)-z_i(t)) \times (z_{i+1}(t)-z_i(t))\big|$.
Thus, we can proceed as in the 2-$D$ case to obtain
\begin{align*}
|f_{rot,3d}(t,x)|~&\leq~ \sum_{i=2}^{N-1} \Xi_i(t)\,|v_{i-1}(t)|\,|z_i(t) - z_{i-1}(t)|\, \big(\xi_{i-1}(t, x)+\xi_i(t, x)\big)\\
&~~~~~~~~+ \sum_{i=2}^{N-1} \Xi_i(t)\,|v_{i-1}(t)| \,{|z_{i-1}(t)-z_i(t)|^2\over |z_{i+1}(t)-z_i(t)|}\,\big(\xi_i(t, x)+\xi_{i+1}(t, x)\big)\\
&\leq~ 4\,  \left(\sum_{i=2}^{N-1} |v_{i-1}(t)|\right) \max_{i=2,\ldots,N}\|z_i-z_{i-1}\|^3_{\con([0,T]\,;\R^d)}\,,
\end{align*}
where we have used the fact that $\Xi_i(t)\leq |z_i(t)-z_{i-1}(t)|\,|z_i(t)-z_{i+1}(t)|$. Therefore, the required inequality follows as in the case of $f_{rot,2d}$. This completes the proof.~~$\diamond$

\subsection{Proof of Lemma~\ref{lem:T}}

The operator ${\bf T}$ is well defined for all times $T^*>0$ satisfying the assumptions of this lemma, thanks to  Lemma~\ref{lem:ode}, and its images  satisfy the properties~\eqref{body2}. Therefore, it remains to prove that it is bounded,
continuous and compact.

\medskip

\n{\bf Step 1: Proof of the continuity of ${\bf T}$.}

\n The proof of this part is similar to that of Lemma~\ref{lem:ode}. Let $u^{(1)},u^{(2)} \in {\cal B}_q$  be fixed. Denote 
$ {\bf T}u^{(j)} = z^{(j)} = (z_1^{(j)}, \ldots, z_N^{(j)}), \, j=1,2$. Then we have for all $t\in [0,T^*]$:
\begin{align*}
|z_i^{(1)}(t)-z_i^{(2)}(t)|& \leq \,{1\over\meas(S(0))}\,\left| \int_0^t \int_{S(z_i^{(1)}(\tau))} \!\!u^{(1)}(\tau,x)\,dx-\int_{S(z_i^{(2)}(\tau))} \!\!u^{(2)}(\tau,x)\,dx\, \,d\tau \right| \\
& \leq \,{1\over\meas(S(0))}\,  \left| \int_0^t \int_\Omega u^{(1)}(\tau,x)-u^{(2)}(\tau,x)  \, \xi_{1,\tau}(x) dx\,d\tau \right| \\
&~~~~~~~~+ \,{1\over\meas(S(0))}\,\left| \int_0^t\int_\Omega u^{(2)}(\tau,x)\,\big( \xi_{1,\tau}(x)-\xi_{2,\tau}(x)\big)\,dx\,\,d\tau \right| \,,
\end{align*}
where we have denoted with $\xi_{j,\tau}(x)$ the characteristic functions of $S(z_i^{(j)}(\tau))$. By applying H\"older inequality to the first term and by repeating for the second term the argument used to derive~\eqref{eq:estimate_with_slices}-(\ref{Tq})  in Lemma~\ref{lem:ode}, we obtain:
\begin{align*}
|z_i^{(1)}(t)-z_i^{(2)}(t)|&\leq \sqrt{{T^*\over\meas(S(0))}}~~\|u^{(1)}-u^{(2)}\|_{\Lsp^2(0,T^*\,;H)}\\
&~~~~~~~~+\,{c_\Omega\,c'_\Omega\,{\cal K}_S\,\sqrt{T^*}\over\meas(S(0))}\,
\|u^{(2)}\|_{\Lsp^2(0,T^*\,;V)}\,\|z_i^{(1)}-z_i^{(2)}\|_{\con([0,T^*]\,;\R^d)}\\
&\leq \,C_\Omega\,\sqrt{{T^*\over\meas(S(0))}}~~\|u^{(1)}-u^{(2)}\|_{\Lsp^2(0,T^*\,;V)}\\
&~~~~~~~~+\,{c_\Omega\,c'_\Omega\,{\cal K}_S\,q\,\sqrt{T^*}\over\meas(S(0))}\,\|z_i^{(1)}-z_i^{(2)}\|_{\con([0,T^*]\,;\R^d)}.
\end{align*}
After taking the supremum over $[0,T^*]$ and  rearranging the terms, we arrive next at the estimate
$$
\left(1-\,{c_\Omega\,c'_\Omega\,{\cal K}_S\,q\,\sqrt{T^*}\over\meas(S(0))}\right)
\|z_i^{(1)}-z_i^{(2)}\|_{\con([0,T^*]\,;\R^d)}
\leq
\,C_\Omega\,\sqrt{{T^*\over\meas(S(0))}}~~\|u^{(1)}-u^{(2)}\|_{\Lsp^2(0,T^*\,;V)},
$$
implying ~\eqref{eq:T_cont} and the continuity of ${\bf T}$, due to (\ref{eq:time_traj_T}).

\medskip

\n{\bf Step 2: Proof of the boundedness of ${\bf T}$}.

\n Let $u$ be fixed in ${\cal B}_q$ and consider any $i=1,\ldots,N$. If we denote with $z_i$ the $i$--th component of ${\bf T}u$, 
we obtain for all $t\in [0,T^*]$ that
\begin{align*}
|z_i(t)|& \leq |z_{i,0}|+ \,{1\over\meas(S(0))}\, \left| \int_0^t \int_{S(z_i(\tau))} \!\! u(\tau,x)\,dx\,d\tau \right| \
&\leq |z_0|+\sqrt{{T^*\over\meas(S(0))}}~~\|u\|_{\Lsp^2(0,T^*\,;H)}.
\end{align*}
Thus,
$$
\|{\bf T}u\|_{[\con([0,T^*]\,;\R^d)]^N}\leq \sqrt{N}\left(|z_0|+\,C_\Omega\,\sqrt{{T^*\over\meas(S(0))}}~~\|u\|_{\Lsp^2(0,T^*\,;V)}\right)\,.
$$

\medskip

\n{\bf Step 3:  Proof of the compactness of ${\bf T}$.} 

\n We will show that any bounded sequence $\{u^{(j)}\}_{j\in\N}$ from $ {\cal B}_q$ is mapped by $ {\bf T}$ to a sequence $\{{z^{(j)} = \bf T} u^{(j)}\}_{j\in\N}$, which has a converging subsequence. To this end, we will prove that $\{{\bf T} u^{(j)}\}_{j\in\N}$ is bounded and equicontinuous. Then Ascoli--Arzel\`a theorem will provide the desirable result.

\n The uniform boundedness of $\{{\bf T} u^{(j)}\}_{j\in\N}$ follows immediately from Step 2. To prove equicontinuity,  fix $t,t+h\in [0,T^*]$ and assume that $h>0$  (the case $h<0$  is analogous). Then, for every index $i=1,\ldots,N$, the $i$--th component $z_i^{(j)}$ of ${\bf T}u^{(j)}$ satisfies:
\begin{align*}
|z_i^{(j)}(t+h)-z_i^{(j)}(t)| &= \,{1\over\meas(S(0))}\,\left|\int_t^{t+h} \int_{S(z_i^{(j)}(\tau))} \!\!u^{(j)}(\tau,x)\,dx\,d\tau\,\right|\\
&\leq\,{C_\Omega\, \|u^{(j)}\|_{\Lsp^2(0,T^*\,;V)}\over \sqrt{\meas(S(0))} }\,\sqrt{h}
\leq\,{ C_\Omega\, q \over \sqrt{\meas(S(0))} }\,\sqrt{h}\,.
\end{align*}
This  completes the proof of Lemma~\ref{lem:T}.~~$\diamond$

\begin{rem}
The results similar to Lemmas ~\ref{lem:ode} and ~\ref{lem:T} were proven in \cite{Kh1}-\cite{KhBook} and \cite{Kh5} for swimming models similar to (\ref{eq:nse})-(\ref{eq:ode})  but for the case of nonstationary Stokes fluid equations. However,  \cite{Kh1}-\cite{KhBook} and \cite{Kh5} deal with  the case when $ u $ is more regular, namely, when $u\in L^2(0,T\,; [L^\infty (\Omega)]^d) \subset  L^2(0,T\,; [H^2(\Omega)]^d)$.
\end{rem}

\subsection{Proof of Lemma~\ref{lem:F}}

\n
The operator ${\bf F}$ is well defined  under the assumptions of Lemma~\ref{lem:F}, thanks to Lemma~\ref{lem:T}.  

\n
In turn, the estimates (\ref{eq:F_bdd})-(\ref{eq:F_cont}) can be derived by straightforward transformations of  explicit elementary expressions in the forcing term (\ref{eq:forces}) (see also  the proof of Lemma~\ref{lem:forces}). Moreover, 
the results very similar to Lemma~\ref{lem:F} were derived in all technical detail in \cite{Kh1}, \cite{KhBook} (Section 12.4) for the 2-$D$ case and in \cite{Kh5} for the 3-$D$ case for swimming models similar to (\ref{eq:nse})-(\ref{eq:ode})  but for the case of nonstationary linear Stokes fluid equations. The difference between the forcing terms in \cite{Kh1}-\cite{KhBook} and \cite{Kh5} and in this paper is that
in the former $ \kappa_i$'s were positive numbers and $ v_i$'s were elements of $ L^\infty (0, T)$, while in this paper 
 $\kappa_i$'s and $ v_i$'s  lie in $ L^2(0,T)$. Furthermore, in the case of 3-$D$ swimming models the formula for rotational forces in  \cite{Kh5}  has a more general  shape but  designed to work only  when the respective triplets $ \{z_{i-1}, z_i, z_{i+1}\}$ are not  in the aligned position (see subsection \ref{sec:force_model} in the above).

\subsection{Proof of Lemma~\ref{lem:S}}

The classical results in Theorems~\ref{thm:existence2d} and~\ref{thm:existence3d} about  solutions to incompressible Navier--Stokes equations~\eqref{eq:nse}, cited in  Section 6 below,  ensure that operator ${\bf S}$ is well defined. 
Let us also remind the reader that  the initial datum $u_0$ is  fixed in our proofs.

\medskip

\n{\it Proof of the boundedness of ${\bf S}$.} 

\n For any  $f\in\Lsp^2(0,T^*\,;  [\Lsp^2(\Omega)]^d)$, estimate~\eqref{eq:S_bdd} follows immediately from~\eqref{eq:l22_norm_2d}  in the 2-$D$  case  and from~\eqref{eq:3dest} in the 3-$D$  case.

\medskip

\n{\it Proof of the continuity of ${\bf S}$.} 

\n {\bf Step 1.} Consider first the case when $d=2$. Let  $(u_1={\bf S}f_1, u_2={\bf S}f_2)$ be  a pair of solutions  to uncoupled (\ref{eq:nse}) generated, given $ u_0$,   by  some  $f_1, f_2 \in \Lsp^2(0,T^*\,; [\Lsp^2(\Omega)]^2)$. To derive \eqref{eq:S_cont2d}, it is sufficient to establish the following estimate:
\begin{equation}\label{eq:estimate2d}
\|u_1-u_2\|^2_{\Lsp^2(0,T^*\,;V)}\leq\,{4\, C_\Omega^2\over \nu^2}\,
e^{{8\over \nu}\,\max\{\|u_1\|^2_{\Lsp^2(0,T^*\,;V)},\|u_2\|^2_{\Lsp^2(0,T^*\,;V)}\}}
\,\|f_1-f_2\|^2_{\Lsp^2(0,T^*\,; [\Lsp^2(\Omega)]^2)}\,.
\end{equation}
where $C_\Omega$ is the constant from Poincar\'e inequality.

\smallskip

\n {\bf Step 2.} From the proof of Theorem~III.3.2 in~\cite{Temam}, we obtain that for any $t\in (0,T^*)$ the difference $u_1(t)-u_2(t)$ satisfies the following equality (see in particular (3.65) in~\cite[Chapter~III]{Temam} or (5), page 145 in \cite{Lad}):
\begin{align*}
{d\over dt}\,\left({1\over 2}\,\|u_1(t)-u_2(t)\|_{\Lsp^2}^2\right)&+\nu\,\|u_1(t)-u_2(t)\|_{H^1_0}^2\\
&=~ \,\langle f_1(t)-f_2(t), u_1(t)-u_2(t)\rangle_{\Lsp^2} \\
&~~~~~~~~ - \, b\big(u_1(t)-u_2(t),u_2(t),u_1(t)-u_2(t)\big)\,,
\end{align*} 
where $\langle \cdot, \cdot \rangle_{\Lsp^2} $ stands for the scalar product in $ \Lsp^2(\Omega)$. Hence, we have
\begin{equation}\label{aux1}
{d\over dt}\,\left(\|u_1(t)-u_2(t)\|_{\Lsp^2}^2\right) +2\nu\,\|u_1(t)-u_2(t)\|_{H^1_0}^2
\end{equation}
$$
\leq~ 2\,\|f_1(t)-f_2(t)\|_{\Lsp^2} ~\|u_1(t)-u_2(t)\|_{\Lsp^2}
$$
$$
~~~~+ 2\, \Big|b\big(u_1(t)-u_2(t),u_2(t),u_1(t)-u_2(t)\big)\Big|
$$
$$
\leq~ 2\,C_\Omega\,\|f_1(t)-f_2(t)\|_{\Lsp^2}~ \|u_1(t)-u_2(t)\|_{H^1_0}
$$
$$
~~~~+ 2\sqrt{2}~\|u_2(t)\|_{H^1_0}~\|u_1(t)-u_2(t)\|_{\Lsp^2}~\|u_1(t)-u_2(t)\|_{H^1_0}\,,
$$
$$
\leq~{ 2\,C_\Omega^2\over \nu}\,\|f_1(t)-f_2(t)\|_{\Lsp^2}^2+\,{\nu\over 2}\, \|u_1(t)-u_2(t)\|_{H^1_0}^2
$$
$$
~~~~+\,{ 4\over \nu}\, \|u_2(t)\|_{H^1_0}^2~\|u_1(t)-u_2(t)\|_{\Lsp^2}^2+\,{\nu\over 2}\, \|u_1(t)-u_2(t)\|_{H^1_0}^2\,,
$$
where $C_\Omega$ is the constant in Poncar\'e inequality and, to evaluate  $|b(\cdot,\cdot,\cdot)|$, we used  (\ref{b2dest}) from Section 6.
Therefore, there holds:
\begin{align}\label{eq:pre_estimate}
{d\over dt}\,\left(\|u_1(t)-u_2(t)\|_{\Lsp^2}^2\right)&+\nu\,\|u_1(t)-u_2(t)\|_{H^1_0}^2\nonumber\\
&\leq \,{2\, C_\Omega^2\over \nu}\,\|f_1(t)-f_2(t)\|_{\Lsp^2}^2+\,{4\over \nu}\,\|u_2(t)\|_{H^1_0}^2~\|u_1(t)-u_2(t)\|_{\Lsp^2}^2\,.
\end{align}
Owing to~\eqref{eq:pre_estimate}, we immediately see that
$$
{d\over dt}\,\left(\|u_1(t)-u_2(t)\|_{\Lsp^2}^2\right)\leq \,{2\, C_\Omega^2\over \nu}\,\|f_1(t)-f_2(t)\|_{\Lsp^2}^2+\,{4\over \nu}\,\|u_2(t)\|_{H^1_0}^2~\|u_1(t)-u_2(t)\|_{\Lsp^2}^2\,.
$$
Hence, an application of Gronwall's inequality yields that
\begin{equation}\label{diff}
\|u_1(t)-u_2(t)\|_{\Lsp^2}^2\leq \,{2\, C_\Omega^2\over \nu}\, \int_0^t \|f_1(s)-f_2(s)\|_{\Lsp^2}^2\,ds\cdot
\exp\left(\,{4\over \nu}\|u_2\|^2_{ \Lsp^2(0,t\,;V)}\right), \;\; t \in [0,T^*].
\end{equation}
In turn, integrating~\eqref{eq:pre_estimate}, we deduce from (\ref{diff}) that
\begin{align*}
\nu\|u_1-u_2\|^2_{\Lsp^2(0,T^*\,;V)}&\leq \|u_1(T^*)-u_2(T^*)\|_{\Lsp^2}^2+\nu\|u_1-u_2\|^2_{\Lsp^2(0,T^*\,;V)}\\
&\leq \,{2\, C_\Omega^2\over \nu}\,\int_0^t \|f_1(s)-f_2(s)\|_{\Lsp^2}^2\,ds+ \max_{s\in[0,T^*]} \|u_1(s)-u_2(s)\|_{\Lsp^2}^2\,\cdot\,{4\over \nu}\, \|u_2\|^2_{\Lsp^2(0,T^*\,;V)}\,\\
&\leq \,{2\, C_\Omega^2\over \nu}\,\int_0^t \|f_1(s)-f_2(s)\|_{\Lsp^2}^2\,ds
\left(1+\,{4\over \nu}\,\|u_2\|^2_{\Lsp^2(0,T^*\,;V)}\,e^{{4\over \nu}\,\|u_2\|^2_{\Lsp^2(0,T^*\,;V)}}\right)\,.
\end{align*}
From here, making use of estimate $  (1+ x e^x) < 2 e^{2x}$ for  $  x \geq 0$, we arrive at
\begin{equation}\label{eq:estimate}
\|u_1-u_2\|^2_{\Lsp^2(0,T^*\,;V)}\leq\,{4\, C_\Omega^2\over \nu^2}\,
e^{{8\over \nu}\,\|u_2\|^2_{\Lsp^2(0,T^*\,;V)}}
\,\|f_1-f_2\|^2_{\Lsp^2(0,T^*\,;[\Lsp^2 (\Omega)]^2))}\,,
\end{equation}
implying~\eqref{eq:S_cont2d}. 
This completes the proof in the $2$d case.

\begin{rem}
An estimate similar to (\ref{eq:estimate2d}) can also be found in \cite{Lad69}, see Theorem 11 and  estimates (45) and (48) on pp. 170-171.
\end{rem}

\n {\bf Step 3: Case $d = 3$.} Similar to Step 2, consider  any two
solutions $u_i={\bf S}f_i, i = 1,2 $ to (\ref{eq:nse}), generated by some  $f_1,f_2$ in $\Lsp^2(0,T^*\,;H)$.
Combining  estimate (11)  from \cite{Lad69}, p. 147, namely:
$$
 \|  \psi  \|_{L^4 (\Omega)}  \leq D^* \|  \psi  \|_{H^1_0}   \;\;  \forall  \psi \in H^1_0 (\Omega), \;  {\rm where} \;  D^* \; {\rm is \;  some \;  positive \; constant},  
 $$ 
with the fact that $ u_2 \in C([0,T^*]; V)$ (see Theorem \ref{thm:existence3d} below) allows us to define 
\begin{equation}\label{ct}
 C_2 = \; \max_{t \in [0, T^*]}  \big|  u_2 (t) \big|_{\Lsp^4} \; \leq \;  3^{1/4} D^* \| u_2 \|_{C([0, T^*]; V)},
\end{equation}
where 
$$ \big|  u_2 (t) \big|_{\Lsp^4} = \bigg[\sum_{i=1}^3 \int_\Omega u_{2i}^4 (t,x)dx \bigg]^{1/4}, \; u_2 = (u_{21}, u_{22}, u_{23}).
$$

\n
To evaluate $ b\big(u_1(t)-u_2(t),u_2(t),u_1(t)-u_2(t)\big)$ in the 3-$D$-case, we will use the chain of  estimates from \cite{Lad69}, p. 145 as follows. 
$$
\Big|b\big(u_1(t)-u_2(t),u_2(t),u_1(t)-u_2(t)\big)\Big|
$$
$$
\leq   \sqrt{3}  \big |  u_2 (t) \big|_{\Lsp^4}  \;  \|  u_1(t)-u_2(t) \|_{H^1_0}  \;   \big| u_1(t)-u_2(t)\big|_{\Lsp^4} 
$$
$$\leq \; \sqrt{3}  C_2  \;  \|  u_1(t)-u_2(t) \|_{H^1_0}  \;    \big| u_1(t)-u_2(t)\big|_{\Lsp^4} . 
$$
Then, making use of estimate (12)  from \cite{Lad69}, p. 147 and (5) on p. 10, namely:
$$
 \|  \psi  \|_{L^4 (\Omega)}  \leq  \delta \|  \psi  \|_{H^1_0 (\Omega)}  + \frac{3^{3/4}}{\delta^3} \|  \psi \|_{L^2 (\Omega)}, \;\;\;\; \psi \in H^1_0 (\Omega), \; \delta> 0,
 $$ 
 we obtain as in  \cite{Lad69}, p. 145:
$$
\Big|b\big(u_1(t)-u_2(t),u_2(t),u_1(t)-u_2(t)\big)\Big|
$$
$$
\leq \;  \sqrt{3} C_2 \|  u_1(t)-u_2(t) \|_{H^1_0}  \bigg[ \varepsilon \|  u_1(t)-u_2(t) \|_{H^1_0} +   C_\varepsilon   \|  u_1(t)-u_2(t)  \|_{[L^2 (\Omega)]^3} \bigg],
$$
$$
\leq \; 2 \sqrt{3} C_2 \varepsilon  \|  u_1(t)-u_2(t) \|^2_{H^1_0}   +   \sqrt{3}   C_2  \frac{C_\varepsilon^2}{4\varepsilon}  \|  u_1(t)-u_2(t)  \|^2_{[L^2 (\Omega)]^3}, \;\; \varepsilon >0,
$$
where  $C_\varepsilon =  3(3)^{3/4}/\varepsilon^3$.

\n
Select now $ \varepsilon = (2 \sqrt{3}    C_2 \} )^{-1} (\nu/4)$, then:
$$
\Big|b\big(u_1(t)-u_2(t),u_2(t),u_1(t)-u_2(t)\big)\Big|
$$
$$
\leq \; \frac{\nu}{4}  \|  u_1(t)-u_2(t) \|^2_{H^1_0} +   \frac{6C_2^2 C_\varepsilon^2}{\nu}    \|  u_1(t)-u_2(t)  \|^2_{[L^2 (\Omega)]^3},
$$
\begin{equation}\label{b3d}
\leq \; \frac{\nu}{4}  \|  u_1(t)-u_2(t) \|^2_{H^1_0} +   D_*    \|  u_1(t)-u_2(t)  \|^2_{[L^2 (\Omega)]^3},
\end{equation}
where we can set, due  to  (\ref{3df}) and (\ref{ct}),
\begin{equation}\label{D*}
D_* \;  =  \hat{C} \hat{L}_S^8 \left( ||  u_0 ||_{V} + || f_2||_{\Lsp^2 (0,T_1\,;[\Lsp^2 (\Omega)]^3)} \right)^8,
\end{equation} 
and  $\hat{C} $ depends on $\nu$. 
 Similar to (\ref{aux1}), we obtain:
 \begin{equation}\label{aux3d}
{d\over dt}\,\left(\|u_1(t)-u_2(t)\|_{\Lsp^2}^2\right) +2\nu\,\|u_1(t)-u_2(t)\|_{H^1_0}^2
\end{equation}
$$
\leq~ 2\,\|f_1(t)-f_2(t)\|_{\Lsp^2} ~\|u_1(t)-u_2(t)\|_{\Lsp^2}
~+ 2\, \Big|b\big(u_1(t)-u_2(t),u_2(t),u_1(t)-u_2(t)\big)\Big|
$$
$$
\leq~{ 2\,C_\Omega^2\over \nu}\,\|f_1(t)-f_2(t)\|_{\Lsp^2}^2+\,{\nu\over 2}\, \|u_1(t)-u_2(t)\|_{H^1_0}^2
$$
$$
~~~~+\,\; \frac{\nu}{2}  \|  (u_1(t)-u_2(t)) \|^2_{H^1_0} +   2 D_*  \|  (u_1(t)-u_2(t)) \|^2_{[L^2 (\Omega)]^3}\,.
$$
Respectively, instead of  (\ref{eq:pre_estimate}), we have:
\begin{align}\label{aux3d2}
{d\over dt}\,\left(\|u_1(t)-u_2(t)\|_{\Lsp^2}^2\right)&+\nu\,\|u_1(t)-u_2(t)\|_{H^1_0}^2\nonumber\\
&\leq \,{2\, C_\Omega^2\over \nu}\,\|f_1(t)-f_2(t)\|_{\Lsp^2}^2+\,2 D_*   \|  (u_1(t)-u_2(t)) \|^2_{[L^2 (\Omega)]^3}\,,
\end{align}
and, furthermore, in place of (\ref{eq:estimate}):
\begin{equation}\label{eq:estimate3d}
\|u_1-u_2\|^2_{\Lsp^2(0,T^*\,;V)}\leq\,{4\, C_\Omega^2\over \nu^2}\,
e^{{4D_*   T^*}}
\,\|f_1-f_2\|^2_{\Lsp^2(0,T^*\,;[\Lsp^2 (\Omega)]^2))}\,,
\end{equation}
implying (\ref{eq:S_cont3d}) in view of (\ref{D*}).
This completes  the proof of Lemma \ref{lem:S}.~~$\diamond$

\section{Appendix A: Solutions to uncoupled Navier--Stokes equations}\label{sec:sol2nse}

In this section we will remind the reader several classical results relevant to the well-posedness of  Navier--Stokes equations which are
used in the proofs of the main results of this paper.

\n
Assume that in \eqref{eq:nse}:
$$
f \in \Lsp^2(0,T\,;[\Lsp^2 (\Omega)]^d), \;\; u_0\in [\Lsp^2(\Omega)]^d.
$$

\medskip

\n
Following~\cite[Chapter~III]{Temam},  we introduce the following trilinear continuous map on $[H^1_0(\Omega)]^d\times [H^1_0(\Omega)]^d\times [H^1_0(\Omega)$ (see Lemma~II.1.1 in~\cite{Temam}, page 161-162 for $ d = 2, 3$, page 159):
\begin{equation}\label{eq:nonlinearity}
b(u,v,w):=\sum_{i,j=1}^d\int_\Omega u_i\,(\D_i v_j)\,w_j\,dx.
\end{equation}
Furthermore, due to Lemma
3.4 (page 292) in \cite[Chapter~III]{Temam}, for $d=2$, for all $u,v,w\in [H^1_0(\Omega)]^2$ there holds
\begin{equation}\label{b2dest}
\big|b(u,v,w)\big|\leq\sqrt{2}~\big(\|u\|_{H^1_0}~\|u\|_{\Lsp^2}\big)^{1\over 2}~\|v\|_{H^1_0}~\big(\|w\|_{H^1_0}~\|w\|_{\Lsp^2}\big)^{1\over 2}\,.
\end{equation}

\medskip

\n
Our results concerning the well--posedness of the swimming model will heavily rely on the following classical results about existence and uniqueness of weak solutions for~\eqref{eq:nse}.

\begin{thm}[Theorems~III.3.1 \&~III.3.2 in~\cite{Temam}, pages 282, 294]\label{thm:existence2d} Let $ T$ be any positive number, $d=2$, $f\in\Lsp^2(0,T\,;[\Lsp^2 (\Omega)]^2)$ and $u_0\in H$.
Then, there exists a unique solution $u\in C([0,T];H)\cap\Lsp^2(0,T\,;V)$ of~\eqref{eq:nse} on $[0,T]$ and the following estimate holds (see   estimates (45) and (48) in Theorem 11 \cite{Lad69}, pp. 170-171 as well as  our estimates (\ref{eq:pre_estimate})-(\ref{eq:estimate}) when $ u_1 = u$ and $u_2, f_2 = 0$):
\begin{equation}\label{eq:l22_norm_2d}
\|u\|_{C (0,T\,; H)} + \|u\|_{\Lsp^2(0,T\,;V)}\leq   \; L_{\bf S} \left(\,\|u_0\|_{\Lsp^2}+ \,\|f\|_{\Lsp^{2}(0,T\,;[\Lsp^2 (\Omega)]^2)} \right)\,
\end{equation}
for some positive constant $L_{\bf S}$.
\end{thm}

\n
Once solution $u$ is found, the following result allows us to retrieve  $p$.
\begin{prop}[Theorem~V.1.7.1 in~\cite{Sohr}, page 295]\label{prop:pressure} Under the assumptions of Theorem \ref{thm:existence2d}, there exists a function ${\cal P}\in\Lsp^{2}(0,T\,;\Lsp^2(\Omega))$ such that 
$
p =\, {\cal P}_t
$
is an associated pressure of $u$ and 
$$
\nabla p = f - u_t +\nu\Delta u-(u\cdot\nabla)u
$$
in the sense of distributions in $(0,T)\times\Omega$.
\end{prop}

\begin{thm}[Theorem 9, \cite{Lad}, p. 203 and Lemma 9, p. 194]\label{thm:existence3d} Let  $ T > 0$ be any given positive number, $d=3$, $f\in\Lsp^2 (0,T\,;[\Lsp^2 (\Omega)]^3)$, $u_0 \in  V$ and $\partial \Omega$ be of class $ C^2$. Then for any $ T_1$ satisfying 
\begin{equation}\label{3dT1}
T_1 \in (0, K), \;\;  K = \min \left\{ T,  C^o \left[ || u_0 ||^2_{ H^1_0} + || f ||^2_{\Lsp^2 (0,T\,; [\Lsp^2 (\Omega)]^3)} \right]^{-2} \right\}
\end{equation}
($ C^o$ depends on $ \Omega$ and $\nu$ only), \eqref{eq:nse} admits a unique solution $(u,p)$ on  $ (0, T_1)$ such that  $ u_t, \Delta  u, \nabla p \in [\Lsp^2 (Q_{T_1})]^3$ and $ u \in C([0,T_1]; V)$. Furthermore,  the following estimates holds (see (\ref{aux3d})-(\ref{eq:estimate3d}) when $ u_1 = u$ and $u_2, f_2 = 0$):
\begin{equation}\label{eq:3dest}
\|u\|_{C([0,T_1]; H)} \; + \;  \|u\|_{\Lsp^2 (0, T_1; V)} \; \leq \; L_S \left( ||  u_0 ||_{\Lsp^2} + || f ||_{\Lsp^2 (0,T_1\,;[\Lsp^2 (\Omega)]^3)} \right)
\end{equation}
for some positive constant $L_S > 0$ and (see \cite{Lad},  Lemma 9, p. 194, (55)):
\begin{equation}\label{3df}
\|u\|_{C([0,T_1]; V)} \; \; \leq \; \hat{L}_S \left( ||  u_0 ||_{H_0^1} + || f ||_{\Lsp^2 (0,T_1\,;[\Lsp^2 (\Omega)]^3)} \right)
\end{equation}
for some positive constant $\hat{L}_S > 0$.
\end{thm}

\n
The cited Theorem 9 in  \cite{Lad}, p. 203 admits  the following modification.
\begin{thm}[Additonal regularity, \cite{Lad69}, Theorem 17, pp. 183-184]\label{thm:addreg}
If  $u_0 \in [H^{2} (\Omega)]^3\bigcap V$, then  in Theorem \ref{thm:existence3d} solution   $ u $ lies in $  [H^{2,1} (Q_{T_1})]^3 \bigcap C([0,T_1]; V)$.   In turn, if  $u_0 \in [H^{2} (\Omega)]^2\bigcap V$, then   $ u \in [H^{2,1} (Q_{T_1})]^2 \bigcap C([0,T_1]; V)$ in Theorem \ref{thm:existence2d}.
\end{thm}

\section{Appendix B: Sections}\label{sec:sobolev_sec}
In the proof of Lemma \ref{lem:ode} in the above we have used the following results (see  assumption ({\bf H2})  for notations). 

\begin{prop}[Proposition~3.105 in~\cite{AFP}]\label{prop:sections}Let $p\in[1,\infty)$, $w \in \Lsp^p(\Omega)$ and $\nu $ be a unit vector in $\R^n$. If we define the generalized directional  derivative of $w$ in the (unit) direction $\nu$ as a distribution $\D_\nu w\in{\cal D}'(\Omega)$ such that 
$$
\int_\Omega \phi\, \D_\nu w \, dx=-\int_\Omega w \,{\partial \phi\over \partial \nu}\,dx
\qquad\qquad
\forall~\phi\in\con_c^\infty(\Omega)\,,
$$
and assume that $\D_\nu w\in \Lsp^p(\Omega)$, then for every $y\in\Omega_\nu$ we have $w_\nu^y\in W^{1,p}(\Omega_\nu^y)$ and $\D w_\nu^y = (\D_\nu w)_\nu^y$, i.e., the differential of the section $w_\nu^y$ is the section of the directional derivative $\D_\nu w$. 
In particular, if $w\in W^{1,p}(\Omega)$, then $w_\nu^y\in W^{1,p}(\Omega_\nu^y)$ for all vectors $\nu, \; | \nu | = 1 $ and all $y\in\Omega_\nu$.
\end{prop}

\n In the proof of Lemma \ref{lem:ode} we use  the following corollary, where  we denote by ${\cal H}^{n-1}$ the $(n-1)$-dimensional Hausdorff measure over Borel sets of $\R^n$.

\begin{cor}\label{cor:section_est} Let $w\in H^1_0(\Omega)$. Then, the following properties  hold.
\begin{description}
\item{\it (i)} There exists a constant $c_\Omega>0$ such that, for all vectors $\nu, \; | \nu | = 1 $ and all $y\in\Omega_\nu$, we have:
\begin{equation}\label{eq:sobolev_sec}
\|w_\nu^y\|_{\Lsp^\infty(\Omega_\nu^y)}~\leq~ c_\Omega~ \|w_\nu^y\|_{H^1_0(\Omega_\nu^y)}\,.
\end{equation}
\item{\it (ii)} There exists a constant $c'_\Omega>0$ such that, for all vectors $\nu, \; | \nu | = 1 $, we have:
\begin{equation}\label{eq:fubini_sec}
\int_{\Omega_\nu}\|w_\nu^y\|_{H^1_0(\Omega_\nu^y)}\,d{\cal H}^{n-1}(y)~\leq~ c'_\Omega~ \|w\|_{H^1_0(\Omega)}\,.
\end{equation}
\end{description}
\end{cor}

\medskip

\n
{\it Proof:} 
Part {\it (i)} is an immediate consequence of Sobolev embedding theorems for bounded sets $\Omega_\nu^y\subseteq\R$, see e.g.~\cite{Brezis}, and of the equivalence between the norm we defined on $H^1_0(\Omega)$ and the usual one on $H^1(\Omega)$. 
Indeed, the Sobolev constants depends only on the measure of the sections $\Omega_\nu^y$, and the latter is bounded by the diameter of $\Omega$, independently on the direction $\nu$ and on the point $y\in\Omega_\nu$.

\n Concerning {\it (ii)}, it is enough to apply H\"older's inequality
and recall Proposition~\ref{prop:sections}
to obtain:
\begin{align*}
\int_{\Omega_\nu}\|w_\nu^y\|_{H^1_0(\Omega_\nu^y)}\,d{\cal H}^{n-1}(y) &\leq
\sqrt{{\cal H}^{n-1}(\Omega_\nu)}~\Bigg[\int_{\Omega_\nu}\left(\int_{\Omega_\nu^y}|\D w_\nu^y(\alpha)|^2\,d\alpha\right)\,d{\cal H}^{n-1}(y)\Bigg]^{1/2}
\\
&=\sqrt{{\cal H}^{n-1}(\Omega_\nu)}~\Bigg[\int_{\Omega_\nu}\left(\int_{\Omega_\nu^y}|(\D_\nu w)^y_\nu(\alpha)|^2\,d\alpha\right)\,d{\cal H}^{n-1}(y)\Bigg]^{1/2}\\
&=\sqrt{{\cal H}^{n-1}(\Omega_\nu)}~\Bigg[\int_\Omega|\D_\nu w(x)|^2\,dx\Bigg]^{1/2}
\end{align*}
Now, notice that the 1st term in the resulting expression, $\sqrt{{\cal H}^{n-1}(\Omega_\nu)}$, can be estimated, independently on the chosen direction $\nu$, by the $(n-1)$--dimensional volume of a ball of radius $R>0$ with $R$ containing  $\Omega$. For the second term we have, in turn: 
$$
\int_\Omega|\D_\nu w(x)|^2\,dx\leq |\nu|^2\,\sum_{i=1}^n\int_\Omega |\D_i w(x)|^2\,dx=\|w\|^2_{H^1_0(\Omega)}\,,
$$
whence the desirable conclusion follows.~~$\diamond$

\end{document}